\documentclass[english]{amsart}
\usepackage[T1]{fontenc}
\usepackage[latin9]{inputenc}
\usepackage{textcomp}
\usepackage{mathrsfs}
\usepackage{amsthm}
\usepackage{amssymb}
\usepackage{stmaryrd}
\usepackage{graphicx}

\makeatletter

\providecommand{\tabularnewline}{\\}

\numberwithin{equation}{section}
\numberwithin{figure}{section}

\usepackage{upgreek}
\newcommand\restrictbar[1]{\raisebox{-1ex}{|}}

\makeatother

\usepackage{babel}
\begin{document}

\title[Actions of Sets and Groups \& Generalized Affine Spaces]{Notes on Actions of Sets and Groups\\
and Generalized Affine Spaces}

\author{Dan Jonsson}
\begin{abstract}
Some well-known and less well-known or new notions related to group
actions are surveyed. Some of these notions are used to generalize
affine spaces. Actions are seen as functions with values in transformation
monoids.
\end{abstract}

\address{University of Gothenburg, SE-405 30 Gothenburg, Sweden}

\email{dan.jonsson@gu.se}

\maketitle

\section{Introduction}

Every group action is an action of a group, but not every action of
a group is a group action; for example, a group may act on a set as
a semigroup. It even makes sense to say that a group acts on a set
as a set; actually, it makes sense to say that a set without any given
algebraic structure acts on another set \cite{key-5}. These notes
aim to clarify such statements, describing both well-known and less
well-known or new notions related to group actions in a systematic
manner. An application of independent interest, touching on differential
geometry, is also presented.

Section 2 puts group actions into context. We first define actions
(of sets) as a general notion and identify several properties that
actions may have. Then actions of groups and some properties that
may characterize such actions are defined, and group actions, monoidal
actions and premonoidal actions are defined in terms of concepts already
defined in Section 2. A few results relating the notions considered,
selected with the application in Section 3 in mind, are also presented.

The application in question concerns generalized affine spaces. Recall
that an affine space can be defined as a set equipped with a regular
group action \cite{key-2}. In Section 3, the notion of affine spaces
is generalized to ``preaffine'' spaces, using the concept of premonoidal
actions defined in Section 2. Section 3 also generalizes affine spaces
to ``multiaffine spaces'', points to a connection between ``semipreaffine''
and multiaffine spaces, and contains a brief discussion of measures
of ``non-affineness'' of particular preaffine or multiaffine spaces,
connecting with the notions of torsion and curvature in differential
geometry. In Appendix A, we consider how to define a generalized affine
space $X$ intrinsically, using a Malcev operation on $X$.

\section{Actions}

\subsection{Actions of sets}

\paragraph*{1}

Let there exist a set $\mathsf{G}$, a set $X$ and a function
\[
\alpha:\mathsf{G}\rightarrow\mathscr{F}_{_{\!}X},\qquad\mathsf{g}\mapsto\alpha\!\left(\mathsf{g}\right),
\]
where $\mathscr{F}_{_{\!}X}$ is the set of all endofunctions on $X$,
that is, the set of all functions $X\rightarrow X$. We call $\alpha$
an \emph{action} of $\mathsf{G}$ on $X$. When $\alpha$ is fixed,
it is convenient to write $\alpha\!\left(\mathsf{g}\right)$ as $\overline{\mathsf{g}}$
and $\alpha\!\left(\mathsf{G}\right)$ as $\overline{\mathsf{G}}$.
We say that $\overline{\mathsf{G}}$ ($\mathsf{G}$) ``acts on''
$X$ (through $\alpha$). 

The idea motivating this construction is that we intend to let $\mathsf{G}$
be equipped with a binary operation $\square$ corresponding to the
natural binary operation $\circ$ on $\mathscr{F}_{_{\!}X}$, namely
function composition defined by $f\!\circ\!g\!\left(x\right)=f\!\left(g\!\left(x\right)\right)$
or, more pedantically, $f\!\circ\!g\!\left(x\right)=$ $f$\raisebox{-.4ex}{|}$_{\mathrm{Im}(g)}$$\left(g\!\left(x\right)\right)$.
The simplest and strongest correspondence between $\square$\linebreak{}
 and $\circ$ is created by letting $\alpha$ be the identity map
on $\mathscr{F}_{_{\!}X}$ and setting $\square=\circ$. More generally,
we can require that $\alpha\!\left(\mathsf{g}\square\mathsf{h}\right)=\alpha\!\left(\mathsf{g}\right)\circ\alpha\!\left(\mathsf{h}\right)$,
but there are also other ways to construct a link between $\square$
and $\circ$.

$\mathscr{F}_{_{\!}X}$ is a monoid under function composition, the
so-called full transformation monoid of $X$; the identity element
of $\mathscr{F}_{_{\!}X}$ is the identity map $\epsilon_{X}$. Typically,
$\alpha\!\left(\mathsf{G}\right)$ is at least a submonoid of $\mathscr{F}_{_{\!}X}$,
often a subgroup of $\mathscr{F}_{_{\!}X}$ whose elements are bijections
$X\rightarrow X$, in particular automorphisms preserving some structure
on $X$.

In this subsection and the next one, some general notions that can
be defined without reference to the binary operation $\square$ will
be introduced. \smallskip{}

\paragraph*{2}

It is clear that the equation $\overline{\mathsf{g}}\!\left(x\right)=y$
has a unique solution $x\in X$ for any $\overline{\mathsf{g}}\in\overline{\mathsf{G}}$
and $y\in X$ if and only if each $\overline{\mathsf{g}}$ is a bijection
$X\rightarrow X$. 

Now suppose that $\overline{\mathsf{g}}$ is a bijection $X\rightarrow X$.
The set $\mathscr{S}_{\!X}$ of all such bijections is a group under
function composition, so there is some $\overline{\mathsf{g}}^{-1}\in\mathscr{S}_{\!X}$
such that $\overline{\mathsf{g}}\circ\overline{\mathsf{g}}^{-1}=\overline{\mathsf{g}}^{-1}\circ\overline{\mathsf{g}}=\epsilon_{X}$.
Hence, there is some $\phi\in\mathscr{F}_{\negthickspace X}$ such
that
\begin{equation}
\overline{\mathsf{g}}\!\left(\phi\!\left(x\right)\right)=\phi\!\left(\overline{\mathsf{g}}\!\left(x\right)\right)=x\qquad\forall x\in X.\label{eq:cancellable-1-1}
\end{equation}
Note that while $\phi=\overline{\mathsf{g}}^{-1}$ is the inverse
in $\mathscr{S}_{\!X}$ of $\overline{\mathsf{g}}$, $\phi$ need
not be the inverse in $\mathsf{\overline{G}}$ of $\overline{\mathsf{g}}$;
in fact, $\overline{\mathsf{G}}$ need not even be a group. On the
other hand, $\phi$ is the semigroup\emph{ }inverse in $\mathscr{F}_{\negthickspace X}$
of $\overline{\mathsf{g}}$; we have $\phi\circ\overline{\mathsf{g}}\circ\phi=\phi$
and $\overline{\mathsf{g}}\circ\phi\circ\overline{\mathsf{g}}=\overline{\mathsf{g}}$.

Finally, if $\overline{\mathsf{g}}\in\overline{\mathsf{G}}$ is a
function such that there is some $\phi\in\mathscr{F}_{\negthickspace X}$
such that (\ref{eq:cancellable-1-1}) holds, then $\overline{\mathsf{g}}\!\left(x\right)=y$
has a solution $x=\phi\!\left(y\right)$ for any $y\in X$ because
$\overline{\mathsf{g}}\!\left(\phi\!\left(y\right)\right)=y$, and
this solution is unique because if $\overline{\mathsf{g}}\!\left(x\right)=y$
then $x=\phi\!\left(\overline{\mathsf{g}}\!\left(x\right)\right)=\phi\!\left(y\right)$.

Thus, the following conditions are equivalent:
\begin{enumerate}
\item For any $\overline{\mathsf{g}}\in\overline{\mathsf{G}}$, there is
some $\phi\in\mathscr{F}_{\negthickspace X}$ such that (\ref{eq:cancellable-1-1})
holds for any $x\in X$.
\item For any $\overline{\mathsf{g}}\in\overline{\mathsf{G}}$ and $y\in X$,
$\overline{\mathsf{g}}\!\left(x\right)=y$ has a unique solution $x\in X$.
\item $\overline{\mathsf{G}}$ is a subset of $\mathscr{S}_{\!X}$, the
set (group) of all bijections $X\rightarrow X$.
\end{enumerate}
We say that $\alpha$ is a \emph{reversible }action of $\mathsf{G}$
on $X$ if and only if one of these conditions holds, and hence all
conditions hold.\smallskip{}

\paragraph*{3}

A\emph{ unital} action (of sets) is an action $\alpha$ of $\mathsf{G}$
on $X$ such that there is some $\mathsf{g}\in\mathsf{G}$ such that
$\overline{\mathsf{g}}\!\left(x\right)=x$ for any $x\in X$, meaning
that $\epsilon_{X}=\overline{\mathsf{g}}\in\overline{\mathsf{G}}$. 

An\emph{ invertible} action (of sets) is an action $\alpha$ of $\mathsf{G}$
on $X$ such that for every $\mathsf{g}\in\mathsf{G}$ there exists
some $\mathsf{h}\in\mathsf{G}$ such that
\begin{equation}
\overline{\mathsf{g}}\!\left(\overline{\mathsf{h}}\!\left(x\right)\right)=\overline{\mathsf{h}}\!\left(\overline{\mathsf{g}}\!\left(x\right)\right)=x\qquad\forall x\in X.\label{eq:invact}
\end{equation}
A comparison of (\ref{eq:invact}) and (\ref{eq:cancellable-1-1})
reveals that an invertible action is reversible. It follows from (\ref{eq:invact})
that for every $\overline{\mathsf{g}}\in\overline{\mathsf{G}}$ there
is some $\overline{\mathsf{g}}^{-1}\in\overline{\mathsf{G}}$, namely
$\overline{\mathsf{h}}$, such that $\overline{\mathsf{g}}^{-1}$
is the inverse in $\mathscr{S}_{X}$ (but not necessarily in $\overline{\mathsf{G}}$)
of $\overline{\mathsf{g}}$. We can show as in § 2 that (\ref{eq:invact})
holds if and only if for any $\mathsf{g}\in\mathsf{G}$ and $y\in X$
there is some $\mathsf{h}\in\mathsf{G}$ such that the equation $\overline{\mathsf{g}}\!\left(x\right)=y$
has the unique solution $x=\overline{\mathsf{h}}\!\left(y\right)$.

A \emph{closed} action (of sets) is an action $\alpha$ of $\mathsf{G}$
on $X$ such that for any $\mathsf{g},\mathsf{h}\in\mathsf{G}$ there
is some $\mathsf{k}\in\mathsf{G}$ such that $\overline{\mathsf{g}}\!\left(\overline{\mathsf{h}}\!\left(x\right)\right)=\overline{\mathsf{k}}\!\left(x\right)$
for all $x\in X$, or equivalently $\overline{\mathsf{g}}\circ\overline{\mathsf{h}}\in\overline{\mathsf{G}}$. 

Note that if $\alpha$ is closed then $\overline{\mathsf{G}}$ is
a semigroup. We have $\epsilon_{X}\circ\overline{\mathsf{g}}=\overline{\mathsf{g}}\circ\epsilon_{X}=\overline{\mathsf{g}}$
for all $\overline{\mathsf{g}}\in\overline{\mathsf{G}}$, so if $\alpha$
is closed and unital so that $\epsilon_{X}\in\overline{\mathsf{G}}$
then $\overline{\mathsf{G}}$ is a monoid, and if $\alpha$ is closed
and invertible then $\overline{\mathsf{G}}$ is a group because $\epsilon_{X}=\overline{\mathsf{g}}\circ\overline{\mathsf{g}}^{-1}=\overline{\mathsf{g}}^{-1}\circ\overline{\mathsf{g}}\in\overline{\mathsf{G}}$
for all $\overline{\mathsf{g}}\in\overline{\mathsf{G}}$, so that
$\epsilon_{X}$ is the identity element in $\overline{\mathsf{G}}$
and $\overline{\mathsf{g}}^{-1}$ is the inverse in $\overline{\mathsf{G}}$
of $\overline{\mathsf{g}}$.

\subsection{Duals of actions}

\paragraph*{4}

A function $f:\mathsf{G}\rightarrow X$ defines two important set
objects: a subset $f\!\left(\mathsf{G}\right)$ of $X$ and a partition
of $\mathsf{G}$ corresponding to the equivalence relation $\sim$
given by
\[
\mathsf{g}\sim\mathsf{h}\quad\Longleftrightarrow\quad f\!\left(\mathsf{g}\right)=f\!\left(\mathsf{h}\right).
\]
$f^{-1}\!\left(x\right)$, where $x\in f\!\left(\mathsf{G}\right)$,
denotes the equivalence class of all $\mathsf{g}\in\mathsf{G}$ such
that $f\!\left(\mathsf{g}\right)=x$. $f$ is surjective if and only
if $f\!\left(\mathsf{G}\right)=X$, and injective if and only if $f^{-1}\!\left(x\right)$
contains exactly one element of $\mathsf{G}$ for each $x\in f\!\left(\mathsf{G}\right)$,
so $f$ is bijective if and only if both these conditions hold.

Let $\mathscr{F}_{\mathsf{G}X}$ be the set of all functions $\mathsf{G}\rightarrow X$.
Any action $\alpha$ of $\mathsf{G}$ on $X$ has a \emph{dual}
\[
\alpha^{*}:X\rightarrow\mathscr{F}_{\mathsf{G}X}
\]
given by $\alpha^{*}\!\left(x\right)\!\left(\mathsf{g}\right):=\alpha\!\left(\mathsf{g}\right)\!\left(x\right)$
for all $x\in X$ and $\mathsf{g}\in\mathsf{G}$. Each $\alpha^{*}\!\left(x\right)$
is thus a function $\mathsf{G}\rightarrow X$; for convenience we
write $\alpha^{*}\!\left(x\right)$ as $\overline{x}$ when $\alpha^{*}$
(or equivalently $\alpha$) is fixed. The subset $\overline{x}\!\left(\mathsf{G}\right)$
of $X$ is called the \emph{orbit} of $x$, and $\overline{x}^{-1}\!\left(y\right)$,
where $x\in X$ and $y\in\overline{x}\!\left(\mathsf{G}\right)$,
is called the \emph{conduit set} from $x$ to $y$. This conduit set
is a subset of $\mathsf{G}$, namely $\left\{ \mathsf{g}\mid\overline{x}\!\left(\mathsf{g}\right)=y\right\} =\left\{ \mathsf{g}\mid\overline{\mathsf{g}}\!\left(x\right)=y\right\} $.

\textsc{\small{}Remark.}{\small{} Orbits and related concepts are
usually defined in the context of group actions. In that context,
it suffices to consider conduit sets of the form $\overline{x}^{-1}\!\left(x\right)$,
each of which corresponds to the }\emph{\small{}stabilizer group}{\small{}
for $x$; this is a subgroup of the group that acts on $X$. Here
we formulate definitions of concepts of this kind in a context with
much less structure. Actually, the situation considered here is analogous
to that in § 2: there we were concerned with equations of the form
$\overline{\mathsf{g}}\left(x\right)=y$, with $x$ unknown; here
we are concerned with equations of the form $\overline{x}\left(\mathsf{g}\right)=y$,
with $\mathsf{g}$ unknown.}{\small \par}

\smallskip{}

\paragraph*{5}

We shall only consider the most important concepts related to the
notions in the preceding paragraph. An action $\alpha$ of $\mathsf{G}$
on $X$ is said to be ($a$) \emph{transitive}, ($b$) \emph{free},
or ($c$) \emph{regular} if and only if $\overline{x}$ is ($a'$)
\emph{surjective}, ($b'$) \emph{injective}, or ($c'$) \emph{bijective}
for each $x\in X$;\emph{ }or equivalently, if and only if ($a''$)
the orbit $\overline{x}\!\left(\mathsf{G}\right)$ of $x$ is $X$
for each $x\in X$, ($b''$) the conduit set $\overline{x}^{-1}\!\left(y\right)$
contains exactly one element of $\mathsf{G}$ for each $x\in X$ and
each $y\in\overline{x}\!\left(\mathsf{G}\right)$, or ($c''$) conditions
($a''$) and ($b''$) hold; or equivalently, if and only if there
is ($a'''$) \emph{at least one}, ($b'''$) \emph{at most one}, or
($c'''$) \emph{exactly one} $\mathsf{g}\in\mathsf{G}$ such that
$\overline{\mathsf{g}}\!\left(x\right)=y$ for any $x,y\in X$. 

Note that if $\alpha$ is regular then $\overline{\mathsf{g}}$, too,
is uniquely determined by the requirement that $\overline{\mathsf{g}}\!\left(x\right)=y$
for given $x,y\in X$.

Assuming that $\alpha$ is regular and denoting the unique $\mathsf{g}\in\mathsf{G}$
such that $\overline{\mathsf{g}}\!\left(x\right)=y$ by $y/x$, we
have
\begin{equation}
\overline{y/x}\!\left(x\right)=y,\label{eq:reg1}
\end{equation}
\begin{equation}
\overline{\mathsf{g}}\!\left(x\right)\!/x=\mathsf{g}.\label{eq:reg2}
\end{equation}

\paragraph*{$5\dagger$}

If $\alpha$ is a free action and $\overline{\mathsf{g}}\!\left(x\right)=\overline{\mathsf{h}}\!\left(x\right)$
for some $x\in X$ and $\mathsf{g},\mathsf{h}\in\mathsf{G}$ then
$\mathsf{g}=\mathsf{h}$. Indeed, if $\overline{\mathsf{g}}\!\left(x\right)=\overline{\mathsf{h}}\!\left(x\right)=y$
for some $x\in X$ but $\mathsf{g}\neq\mathsf{h}$ then $\alpha$
is not a free action. Hence, if $\alpha$ is a free action then the
following conditions are equivalent, since obviously (2) implies (3)
and (3) implies (1).
\begin{enumerate}
\item $\overline{\mathsf{g}}\!\left(x\right)=\overline{\mathsf{h}}\!\left(x\right)$
for some $x\in X$.
\item $\mathsf{g}=\mathsf{h}$.
\item $\overline{\mathsf{g}}=\overline{\mathsf{h}}$.
\end{enumerate}
In particular, $\overline{\mathsf{g}}=\overline{\mathsf{h}}$ implies
$\mathsf{g}=\mathsf{h}$, so a free action is an injective function.\pagebreak{}

\subsection{Actions of groups}

\paragraph*{6}

In this subsection, we let the set $\mathsf{G}$ be a group $G$ with
generic elements denoted $g,h,\ldots$ and identity element $e$.
The binary operation $\square$ in $G$ linked to function composition
$\circ$ in $\mathscr{F}_{X}$ by an action $\alpha$ is simply group
multiplication.

A \emph{unital action of groups} is an action $\alpha$ of $G$ on
$X$ such that 
\[
\overline{e}\!\left(x\right)=x\qquad\forall x\in X,
\]
or $\overline{e}=\epsilon{}_{\!X}$. A unital action of groups is
trivially unital as an action of sets.

An \emph{invertible action of groups} is an action $\alpha$ of $G$
on $X$ such that for every $g\in G$
\begin{equation}
\overline{g}\!\left(\overline{g^{-1}}\!\left(x\right)\right)=\overline{g^{-1}}\!\left(\overline{g}\!\left(x\right)\right)=x\qquad\forall x\in X.\label{eq:covact}
\end{equation}
A comparison of (\ref{eq:covact}) and (\ref{eq:invact}) shows that
an invertible action of groups is invertible as an action of sets,
and hence a reversible action. Thus, $\overline{g^{-1}}$ is the inverse
in $\mathscr{S}_{\!X}$ of $\overline{g}$ (§ 2), but not necessarily
the inverse in $\overline{G}$ of $\overline{g}$ even though $\overline{g^{-1}}\in\overline{G}$.
(There are trivial counterexamples for $G=\left\{ e\right\} $.) Again
reasoning as in § 2, we see that $\alpha$ is an invertible action
of groups if and only if the equation $\overline{g}\!\left(x\right)=y$
has the unique solution $x=\overline{g^{-1}}\!\left(y\right)$ for
every $g\in G$ and $y\in X$. 

A \emph{closed action of groups} (or \emph{semigroup action} of groups)
is an action of $G$ on $X$ such that 
\[
\overline{g}\!\left(\overline{h}\!\left(x\right)\right)=\overline{gh}\!\left(x\right)\qquad\forall g,h\in G,\forall x\in X,
\]
or $\overline{g}\circ\overline{h}=\overline{gh}$. A closed action
of groups is clearly closed as an action of sets. 

It is obvious what to mean by transitive, free or regular actions
of groups for various kinds of actions of groups. For example, a regular
unital action of groups is a unital action of groups that is regular
as an action of sets.\smallskip{}

\paragraph*{$6\dagger$}

If $\alpha$ is a free unital action of groups and $\overline{g}\!\left(x_{0}\right)=x_{0}$
for some $g\in G$ and $x_{0}\in X$, so that $\overline{g}\!\left(x_{0}\right)=\epsilon_{\!X}\!\left(x_{0}\right)=\overline{e}\!\left(x_{0}\right)$,
then $g=e$ by § $5\dagger$, and $\overline{g}=\overline{e}=\epsilon_{X}$.

\subsection{Monoidal and premonoidal actions of groups; translations}

\paragraph*{7}

If $\alpha$ is a closed action of groups then $\overline{G}$ is
a group. To verify this, first recall that $\overline{g}\circ\overline{h}=\overline{gh}\in\overline{G}$
for any $\overline{g},\overline{h}\in\overline{G}$, so $\overline{G}$
is a semigroup. Also, $\overline{e}\circ\overline{g}=\overline{eg}=\overline{g}=\overline{ge}=\overline{g}\circ\overline{e}$
for all $\overline{g}\in\overline{G}$, so $\overline{e}$ is the
unique identity element in $\overline{G}$. Finally, $\overline{g}\circ\overline{g^{-1}}=\overline{gg^{-1}}=\overline{e}=\overline{g^{-1}g}=\overline{g^{-1}}\circ\overline{g}$
for all $\overline{g}\in\overline{G}$, so $\overline{g^{-1}}$ is
the unique inverse in $\overline{G}$ of $\overline{g}$. Thus, $\alpha$
is a group homomorphism $G\rightarrow\overline{G}$, but $\alpha$
need not be a group homomorphism $G\rightarrow\mathscr{S}_{_{\!}X}$.
In particular, the identity element $\overline{e}$ in $\overline{G}$
is not necessarily equal to the identity element $\epsilon{}_{\!X}$
in $\mathscr{F}_{_{\!}X}$ and $\mathscr{S}_{_{\!}X}$. 

\smallskip{}

\textsc{\small{}Example.}{\small{} Consider a set $X=\left\{ a,b,c,d\right\} $
and the functions
\begin{gather*}
\varepsilon:a\mapsto a,\;b\mapsto a,\;c\mapsto c,\;d\mapsto d;\\
\phi:a\mapsto a,\;b\mapsto a,\;c\mapsto d,\;d\mapsto c.
\end{gather*}
 One easily verifies that $\varepsilon\circ\varepsilon=\phi\circ\phi=\varepsilon,\:\varepsilon\circ\phi=\phi\circ\varepsilon=\phi$.
Thus, $\Gamma=\left\{ \varepsilon,\phi\right\} $ is a group under
function composition with identity element $\varepsilon$. Consider
an abstract group $G=\left\{ e,f\right\} $ where $ee=ff=e,\:ef=fe=f$,
so that $e$ is the identity element in $G$. Let
\[
\alpha:e\mapsto\varepsilon,\;f\mapsto\phi
\]
be an action of $G$ on $X$. Then $\alpha$ is a closed action of
groups, $\alpha\!\left(G\right)=\Gamma$ is a group, and $\alpha\!\left(e\right)$
is the identity element $\varepsilon$ in $\Gamma$. Note that $\varepsilon$
is not a bijection. Thus, the identity element in $\Gamma$ is not
the identity element in the monoid $\mathscr{F}_{_{\!}X}$, and also
not the identity element in the group $\mathscr{S}_{_{\!}X}$.\smallskip{}
}{\small \par}

\paragraph*{8}

The previous paragraph suggests that we should look for a notion that
is stronger than that of a closed action of groups. A \emph{monoidal
action} $\alpha$ of $G$ on $X$ is a \emph{unital}, closed action
of groups. Recall that $\overline{G}$ is a group with identity element
$\overline{e}$ and that $\overline{g}\circ\overline{g^{-1}}=\overline{e}=\overline{g^{-1}}\circ\overline{g}$
for all $g\in G$ (§~7). Hence, $\overline{e}=\epsilon{}_{\!X}$
implies that $\alpha$ is an invertible action of groups and therefore
a reversible action of $G$ on $X$, and that $\overline{g^{-1}}$
is the inverse in $\mathscr{S}_{\negthickspace X}$ of $\overline{g}$.
Thus, we have $\alpha\!\left(e\right)=\epsilon_{X}$, $\alpha\!\left(g^{-1}\right)=\alpha\!\left(g\right)^{-1}$,
$\alpha\!\left(gh\right)=\alpha\!\left(g\right)\circ\alpha\!\left(h\right)$,
and $\overline{G}\subset\mathscr{S}_{\negthickspace X}$, so $\alpha$
is a group homomorphism $G\rightarrow\mathscr{S}_{X}$, also known
as a \emph{group action} of $G$ on $X$.

It turns out that a notion that is somewhat weaker than the one just
defined is surprisingly useful, namely a \emph{premonoidal action}
of $G$ on $X$, defined as as a unital action of groups that is closed
as an action of \emph{sets}.

A \emph{regular (pre)monoidal action} is, of course, a (pre)monoidal
action that is regular as an action of sets. \smallskip{}

\paragraph*{9}

Let $\rho$ and $\sigma$ be, respectively, a monoidal and a premonoidal
action of $G$ on $X$. By § 8, $\rho\!\left(G\right)$ is a group
of transformations of $X$ with $\circ$ as binary operation. Although
$\sigma\!\left(G\right)$ need not be a group, a connection between
$G$ and $\sigma\!\left(G\right)$ similar to that between $G$ and
$\rho\!\left(G\right)$ exists if $\sigma$ is a \emph{regular }premonoidal
action of $G$ on $X$.

To show this, first use the fact that $\sigma$ is closed as an action
of sets, which means that for any $g,h\in G$ we have $\overline{g}\circ\overline{h}\in\sigma\!\left(G\right)$,
so $\sigma\!\left(G\right)$ is a semigroup under function composition.
Also, $\sigma$ is a unital action of groups, so $\epsilon_{\!X}=\overline{e}\in\sigma\!\left(G\right)$,
and as $\epsilon_{\!X}\circ\overline{g}=\overline{g}=\overline{g}\circ\epsilon_{\!X}$
for any $\overline{g}\in\sigma\!\left(G\right)$, $\epsilon_{\!X}$
is the identity element in $\sigma\!\left(G\right)$. Finally, $\sigma$
is regular, so for any $\overline{g}\in\sigma\!\left(G\right)$ and
any $x_{0}\in X$ there is a unique $\overline{h}\in\sigma\!\left(G\right)$
such that $\overline{h}\!\left(\overline{g}\!\left(x_{0}\right)\right)=x_{0}$,
and hence $\overline{g}\!\left(\overline{h}\!\left(\overline{g}\!\left(x_{0}\right)\right)\right)=\overline{g}\!\left(x_{0}\right)$.
By assumption, $\overline{h}\circ\overline{g},\,\overline{g}\circ\overline{h}\in\sigma\!\left(G\right)$,
so by § $6\dagger$ we have $\overline{h}\circ\overline{g}=\overline{g}\circ\overline{h}=\epsilon_{\!X}$
since $\sigma$ is a free unital action of groups. Hence, $\overline{h}\in\sigma\!\left(G\right)$
is the inverse in $\sigma\!\left(G\right)$ of $\overline{g}$. We
have thus shown that $\sigma\!\left(G\right)$ is a group of functions
$X\rightarrow X$ under function composition.

Note that $\sigma$ is a reversible action since it is invertible
as an action of sets. We call the bijections in\emph{ $\sigma\!\left(G\right)$}
\emph{translations}, and $\sigma\!\left(G\right)$ is accordingly
called a \emph{translation group} on $X$. As $\sigma$ is a free
action, translations are characterized by the property that $\epsilon_{\!X}$
is the only translation such that $\epsilon_{\!X}\!\left(x_{0}\right)=x_{0}$
for some $x_{0}\in X$; any translation that has a fixed point has
only fixed points. 

Also note that $\sigma$ is injective, since it is free (§ $5\dagger$),
so $\sigma$ is a bijection between $G$ and $\sigma\!\left(G\right)$.
If $\sigma$ happens to be a regular \emph{monoidal} action, or equivalently
a regular group action, then $\sigma$ is an isomorphism between $G$
and $\sigma\!\left(G\right)$, but in general the bijection between
$G$ and $\sigma\!\left(G\right)$ is not an isomorphism. This fact
is illustrated by the following example.

\smallskip{}

\textsc{\small{}Example.}{\small{} Cayley tables describing the cyclic
group of order eight $C_{8}$ and the quaternion group $Q$ are shown
below.\medskip{}
}{\small \par}

\begin{center}
{\small{}}%
\begin{tabular}{c|c|c|c|c|c|c|c|c|}
\multicolumn{1}{c}{{\small{}$C_{8}$}} & \multicolumn{1}{c}{{\small{}$e$}} & \multicolumn{1}{c}{{\small{}$a$}} & \multicolumn{1}{c}{{\small{}$a^{2}$}} & \multicolumn{1}{c}{{\small{}$a^{3}$}} & \multicolumn{1}{c}{{\small{}$a^{4}$}} & \multicolumn{1}{c}{{\small{}$a^{5}$}} & \multicolumn{1}{c}{{\small{}$a^{6}$}} & \multicolumn{1}{c}{{\small{}$a^{7}$}}\tabularnewline
\cline{2-9} 
{\small{}$e$} & {\small{}$e$} & {\small{}$a$} & {\small{}$a^{2}$} & {\small{}$a^{3}$} & {\small{}$a^{4}$} & {\small{}$a^{5}$} & {\small{}$a^{6}$} & {\small{}$a^{7}$}\tabularnewline
\cline{2-9} 
{\small{}$a$} & {\small{}$a$} & {\small{}$a^{2}$} & {\small{}$a^{3}$} & {\small{}$a^{4}$} & {\small{}$a^{5}$} & {\small{}$a^{6}$} & {\small{}$a^{7}$} & {\small{}$e$}\tabularnewline
\cline{2-9} 
{\small{}$a^{2}$} & {\small{}$a^{2}$} & {\small{}$a^{3}$} & {\small{}$a^{4}$} & {\small{}$a^{5}$} & {\small{}$a^{6}$} & {\small{}$a^{7}$} & {\small{}$e$} & {\small{}$a$}\tabularnewline
\cline{2-9} 
{\small{}$a^{3}$} & {\small{}$a^{3}$} & {\small{}$a^{4}$} & {\small{}$a^{5}$} & {\small{}$a^{6}$} & {\small{}$a^{7}$} & {\small{}$e$} & {\small{}$a$} & {\small{}$a^{2}$}\tabularnewline
\cline{2-9} 
{\small{}$a^{4}$} & {\small{}$a^{4}$} & {\small{}$a^{5}$} & {\small{}$a^{6}$} & {\small{}$a^{7}$} & {\small{}$e$} & {\small{}$a$} & {\small{}$a^{2}$} & {\small{}$a^{3}$}\tabularnewline
\cline{2-9} 
{\small{}$a^{5}$} & {\small{}$a^{5}$} & {\small{}$a^{6}$} & {\small{}$a^{7}$} & {\small{}$e$} & {\small{}$a$} & {\small{}$a^{2}$} & {\small{}$a^{3}$} & {\small{}$a^{4}$}\tabularnewline
\cline{2-9} 
{\small{}$a^{6}$} & {\small{}$a^{6}$} & {\small{}$a^{7}$} & {\small{}$e$} & {\small{}$a$} & {\small{}$a^{2}$} & {\small{}$a^{3}$} & {\small{}$a^{4}$} & {\small{}$a^{5}$}\tabularnewline
\cline{2-9} 
{\small{}$a^{7}$} & {\small{}$a^{7}$} & {\small{}$e$} & {\small{}$a$} & {\small{}$a^{2}$} & {\small{}$a^{3}$} & {\small{}$a^{4}$} & {\small{}$a^{5}$} & {\small{}$a^{6}$}\tabularnewline
\cline{2-9} 
\end{tabular}
\par\end{center}{\small \par}

{\small{}\smallskip{}
}{\small \par}

\begin{center}
{\small{}}%
\begin{tabular}{c|c|c|c|c|c|c|c|c|}
\multicolumn{1}{c}{{\small{}$Q$}} & \multicolumn{1}{c}{{\small{}$1$}} & \multicolumn{1}{c}{{\small{}$i$}} & \multicolumn{1}{c}{{\small{}$j$ }} & \multicolumn{1}{c}{{\small{}$k$}} & \multicolumn{1}{c}{{\small{}$-1$}} & \multicolumn{1}{c}{{\small{}$-i$}} & \multicolumn{1}{c}{{\small{}$-j$}} & \multicolumn{1}{c}{{\small{}$-k$}}\tabularnewline
\cline{2-9} 
{\small{}$1$} & {\small{}$1$} & {\small{}$i$} & {\small{}$j$ } & {\small{}$k$} & {\small{}$-1$} & {\small{}$-i$} & {\small{}$-j$} & {\small{}$-k$}\tabularnewline
\cline{2-9} 
{\small{}$i$} & {\small{}$i$} & {\small{}$-1$} & {\small{}$k$} & {\small{}$-j$} & {\small{}$-i$} & {\small{}$1$} & {\small{}$-k$} & {\small{}$j$}\tabularnewline
\cline{2-9} 
{\small{}$j$} & {\small{}$j$} & {\small{}$-k$} & {\small{}$-1$} & {\small{}$i$} & {\small{}$-j$} & {\small{}$k$} & {\small{}$1$} & {\small{}$-i$}\tabularnewline
\cline{2-9} 
{\small{}$k$} & {\small{}$k$} & {\small{}$j$} & {\small{}$-i$} & {\small{}$-1$} & {\small{}$-k$} & {\small{}$-j$} & {\small{}$i$} & {\small{}$1$}\tabularnewline
\cline{2-9} 
{\small{}$-1$} & {\small{}$-1$} & {\small{}$-i$} & {\small{}$-j$} & {\small{}$-k$} & {\small{}1} & {\small{}$i$} & {\small{}$j$} & {\small{}$k$}\tabularnewline
\cline{2-9} 
{\small{}$-i$} & {\small{}$-i$} & {\small{}$1$} & {\small{}$-k$} & {\small{}$j$} & {\small{}$i$} & {\small{}$-1$} & {\small{}$k$} & {\small{}$-j$}\tabularnewline
\cline{2-9} 
{\small{}$-j$} & {\small{}$-j$} & {\small{}$k$} & {\small{}$1$} & {\small{}$-i$} & {\small{}$j$} & {\small{}$-k$} & {\small{}$-1$} & {\small{}$i$}\tabularnewline
\cline{2-9} 
{\small{}$-k$} & {\small{}$-k$} & {\small{}$-j$} & {\small{}$i$} & {\small{}$1$} & {\small{}$k$} & {\small{}$j$} & {\small{}$-i$} & {\small{}$-1$}\tabularnewline
\cline{2-9} 
\end{tabular}{\small{}\medskip{}
}
\par\end{center}{\small \par}

{\small{}Let $\mathscr{T}_{Q}=\left\{ \widehat{1},\;\widehat{i},\;\widehat{j},\;\widehat{k},\;\widehat{-1},\;\widehat{-i},\;\widehat{-j},\;\widehat{-k}\right\} $
be bijections $Q\rightarrow Q$ defined by 
\[
\widehat{z}\!\left(x\right)=zx\qquad\forall x,z\in Q.
\]
$\tau:z\!\mapsto\!\left(x\!\mapsto\!zx\right)$ is a bijection $Q\rightarrow\mathscr{T}_{Q}$
because $\tau$ is clearly surjective and also injective since $zx=z'x$
implies $z=z'$. Furthermore, $\widehat{zy}\!\left(x\right)=\widehat{z}\circ\widehat{y}\!\left(x\right)$
for all $x,y,z\in Q$ because $Q$ is associative, and $\widehat{1}=\epsilon_{Q}$.
$\mathscr{T}_{Q}$ is clearly a group isomorphic to $Q$. }{\small \par}

{\small{}Now consider the function
\begin{gather*}
\sigma:C_{8}\rightarrow\mathscr{F}_{_{\!}Q},\\
e\mapsto\widehat{1},\;a\mapsto\widehat{i},\;a^{2}\mapsto\widehat{j},\;a^{3}\mapsto\widehat{k},\;a^{4}\mapsto\widehat{-1},\;a^{5}\mapsto\widehat{-k},\;a^{6}\mapsto\widehat{-j},\;a^{7}\mapsto\widehat{-i}.
\end{gather*}
$\sigma$ is an action of groups of $C_{8}$ on $Q$, }\emph{\small{}unital}{\small{}
since $\sigma\!\left(e\right)=\widehat{1}=\epsilon_{Q}$; }\emph{\small{}closed}{\small{}
as an action of sets since $\sigma\!\left(a^{m}\right)\circ\sigma\!\left(a^{n}\right)=\widehat{x}\circ\widehat{y}=\widehat{xy}\in\mathscr{T}_{Q}=\sigma\!\left(C_{8}\right)$;
and }\emph{\small{}regular}{\small{} as an action of sets since the
equation $zx=y$ has a unique solution $z=yx^{-1}\in Q$ for any $x,y\in Q$
so that the equation $\sigma\!\left(a^{n}\right)\!\left(x\right)=y$
has a unique solution $a^{n}=\sigma^{-1}\!\left(\widehat{yx^{-1}}\right)\in C_{8}$
for any $x,y\in Q$. Thus, $\sigma$ is a regular premonoidal action
of $C_{8}$ on $Q$, $\sigma$ is a bijection $C_{8}\rightarrow\mathscr{T}_{Q}$
by its definition, and $\mathscr{T}_{Q}$ is a group, Yet, $C_{8}$
and $\mathscr{T}_{Q}$ are not isomorphic, even though in this example,
by design, $\sigma\!\left(g^{-1}\right)=\sigma\!\left(g\right)^{-1}$
for all $g\in C_{8}$. In fact, $C_{8}$ is abelian while $\mathscr{T}_{Q}$,
being isomorphic to $Q$, is non-abelian.}{\small \par}

\subsection{Actions and binary functions}

\paragraph*{10}

For any action $\alpha$ of $\mathsf{G}$ on $X$ there is a unique
corresponding binary function
\begin{equation}
\beta:\mathsf{G}\times X\rightarrow X,\qquad\left(\mathsf{g},x\right)\mapsto\mathsf{g}x,\label{eq:beta}
\end{equation}
alternatively denoted $\widehat{\alpha}$, defined by 
\begin{equation}
\mathsf{g}x:=\overline{\mathsf{g}}\!\left(x\right)\qquad\forall\mathsf{g}\in\mathsf{G},\forall x\in X.\label{eq:g-gbar}
\end{equation}
Conversely, for any binary function $\beta$ of the form (\ref{eq:beta})
there is a unique action $\alpha$, alternatively denoted $\widehat{\beta}$,
defined by 
\begin{equation}
\overline{\mathsf{g}}:=\mathsf{g},\quad\overline{\mathsf{g}}\!\left(x\right):=\mathsf{g}x\qquad\forall\mathsf{g}\in\mathsf{G},\forall x\in X.\label{eq:gbar-g}
\end{equation}

We can thus write $\overline{\mathsf{g}}\!\left(x\right)$ as $\mathsf{g}x$,
and we have the following correspondences:\smallskip{}

\begin{center}
\begin{tabular}{cc}
\noalign{\vskip\doublerulesep}
$\alpha$ & $\widehat{\alpha}$\tabularnewline[\doublerulesep]
\hline 
\noalign{\vskip\doublerulesep}
\noalign{\vskip\doublerulesep}
$\overline{y/x}\!\left(x\right)=y$ & $\left(y/x\right)x=y$\tabularnewline[\doublerulesep]
\noalign{\vskip\doublerulesep}
\noalign{\vskip\doublerulesep}
$\overline{\mathsf{g}}\!\left(x\right)\!/x=\mathsf{g}$ & $\mathsf{g}x/x=\mathsf{g}$\tabularnewline[\doublerulesep]
\noalign{\vskip\doublerulesep}
\noalign{\vskip\doublerulesep}
$\mathsf{\overline{g}}\!\left(x\right)=x$ & $\mathsf{g}x=x$\tabularnewline[\doublerulesep]
\noalign{\vskip\doublerulesep}
\noalign{\vskip\doublerulesep}
$\overline{\mathsf{g}}\!\left(\overline{\mathsf{h}}\!\left(x\right)\right)=\overline{\mathsf{h}}\!\left(\overline{\mathsf{g}}\!\left(x\right)\right)=x$ & $\mathsf{g}\!\left(\mathsf{h}x\right)=\mathsf{h}\!\left(\mathsf{g}x\right)=x$\tabularnewline[\doublerulesep]
\noalign{\vskip\doublerulesep}
\noalign{\vskip\doublerulesep}
$\overline{\mathsf{g}}\circ\overline{\mathsf{h}}\!\left(x\right)=\overline{\mathsf{g}}\!\left(\overline{\mathsf{h}}\!\left(x\right)\right)=\overline{\mathsf{k}}\!\left(x\right)$ & $\mathsf{g}\!\left(\mathsf{h}x\right)=\mathsf{k}x$\tabularnewline[\doublerulesep]
\noalign{\vskip\doublerulesep}
\end{tabular}
\par\end{center}

\smallskip{}

The formal properties of a binary function $\beta$ are similar to
the formal properties of the corresponding action $\widehat{\beta}$.
For example, $\beta$ represents a reversible action $\widehat{\beta}$
if and only if the equation $\mathsf{g}x=y$ has a unique solution
$x$ for any $\mathsf{g}\in\mathsf{G}$ and $y\in X$.

Note that the set of functions $\mathsf{G}\times X\rightarrow X$
is in one-to-one correspondence with the set of functions $\mathsf{G}\rightarrow\mathscr{F}_{\negthickspace X}$
rather than the set of functions $\mathsf{G}\rightarrow\mathscr{S}_{\negthickspace X}$.
This suggests that an action should indeed be defined as a function
$\alpha:\mathsf{G}\rightarrow\mathscr{F}_{\negthickspace X}$ rather
than a function $\alpha:\mathsf{G}\rightarrow\mathscr{S}_{\negthickspace X}$
(if actions are not defined as binary functions).

\smallskip{}

\paragraph*{11}

Let $\mathsf{G}$ be a group $G$ with generic elements $g,h,\ldots$
and identity element $e$, so that we can write (\ref{eq:beta}) as
\begin{equation}
\beta:G\times X\rightarrow X,\qquad\left(g,x\right)\mapsto gx.\label{eq:g-gbar-1}
\end{equation}
The correspondence between $\overline{\mathsf{g}}\!\left(x\right)$
and $\mathsf{g}x$ is now rendered as a correspondence \linebreak{}
between $\overline{g}\!\left(x\right)$ and $gx$, and we can rewrite
the displayed correspondences in the previous paragraph accordingly.
We also have some further correspondences:\smallskip{}

\begin{center}
\begin{tabular}{cc}
\noalign{\vskip\doublerulesep}
$\alpha$ & $\widehat{\alpha}$ \tabularnewline[\doublerulesep]
\hline 
\noalign{\vskip\doublerulesep}
\noalign{\vskip\doublerulesep}
$\overline{e}\!\left(x\right)=x$ & $ex=x$\tabularnewline[\doublerulesep]
\noalign{\vskip\doublerulesep}
\noalign{\vskip\doublerulesep}
$\overline{g}\!\left(\overline{g^{-1}}\!\left(x\right)\right)=\overline{g^{-1}}\!\left(\overline{g}\!\left(x\right)\right)=x\qquad$ & $\qquad g\left(g^{-1}x\right)=g^{-1}\left(gx\right)=x$\tabularnewline[\doublerulesep]
\noalign{\vskip\doublerulesep}
\noalign{\vskip\doublerulesep}
$\overline{g}\circ\overline{h}\!\left(x\right)=\overline{g}\!\left(\overline{h}\!\left(x\right)\right)=\overline{gh}\left(x\right)$ & $g\left(hx\right)=\left(gh\right)x$\tabularnewline[\doublerulesep]
\noalign{\vskip\doublerulesep}
\end{tabular}\smallskip{}

\par\end{center}

In the literature, a group action is mostly defined as a binary function
$\beta$ such that $ex=x$ and $g\left(hx\right)=\left(gh\right)x$
for all $g,h\in G$ and all $x\in X$. A binary function of this kind
represents a monoidal action $G\rightarrow\mathscr{F}_{\negthickspace X}$,
where $\overline{e}\!\left(x\right)=x$ and $\overline{g}\!\left(\overline{h}\!\left(x\right)\right)=\overline{gh}\!\left(x\right)$,
and such an action is invertible and hence reversible (§ 8). Alternatively,
$g\left(g^{-1}x\right)=\left(gg^{-1}\right)x=ex=x$ and similarly
$g^{-1}\left(gx\right)=x$ for all $g\in G$ and $x\in X$, so $gx=y$
has the unique solution $x=g^{-1}y$ for any $g\in G$ and $y\in X$,
implying that $\widehat{\beta}$ is reversible. The assumption $ex=x$
is thus introduced to ensure that $\beta$ represents a group action
$\widehat{\beta}:G\rightarrow\mathscr{S}_{\negthickspace X}$.

\smallskip{}

\paragraph*{{\small{}12}}

{\small{}Let us briefly consider the consequences of defining an action
as a binary function $G\times X\rightarrow X$ rather than a function
$G\rightarrow\mathscr{F}_{\negthickspace X}$, so that $\overline{g}\!\left(x\right)=y$
is replaced by $gx=y$. This definition and notation encourages us
to identify the element $g$ of the group $G$ with the corresponding
function $\overline{g}:X\rightarrow X$ that sends $x$ to $y$, so
we say that $g$ ``acts on'' $x$ or that $G$ ``acts on'' $X$.
This identification tends to work well as long as we are dealing with
a group action, because the identity $\left(gh\right)x=g\left(hx\right)$
means that we can think about the group product $gh$ as the ``product''
$g\circ h$ of $g$ and $h$ regarded as functions. (For contravariant
actions, the ``product'' of $g$ and $h$ is, analogously, $h\circ g$;
see § 14.) This way of thinking may work well even if $\alpha$ is
not injective, because we do not automatically assume that $g_{1}x=g_{2}x$
implies $g_{1}=g_{2}$. }{\small \par}

{\small{}Identifying $g$ and $\overline{g}$ can lead to confusion
when we do not have $\overline{gh}=\overline{g}\circ\overline{h}$
or $\overline{gh}=\overline{h}\circ\overline{g}$, however. For example,
the argument in § 9 depends on the distinction between $g$ and $\overline{g}$
-- recall, in particular, that $C_{8}$ is an abelian group while
$\overline{C_{8}}=\mathscr{T}_{Q}$ is a non-abelian group. In the
application in Section 3, the distinction between $g$ and $\overline{g}$
takes the form of a distinction between the }\emph{\small{}vector}{\small{}
$\mathbf{v}$ (an algebraic concept) and the corresponding }\emph{\small{}translation}{\small{}
$\overline{\mathbf{v}}$ (a geometric concept). In this case, $G$
is the abelian additive group of a vector space, while $\overline{G}$
is a possibly non-abelian group of translations. (In particular, compare
(\ref{eq:aChaslesV}), valid in an affine space, with (\ref{eq:closedset}),
valid in a preaffine space.)}{\small \par}

{\small{}If we define actions as done in § 1, we can still use binary
functions $\mathsf{G}\times X\rightarrow X$ to }\emph{\small{}represent}{\small{}
actions, but this is not necessary, although it may be helpful in
some cases.}{\small \par}

\subsection{On notation}

\paragraph*{13}

The expression $\overline{g}\!\left(x\right)$ uses \emph{left-handed
notation} for actions, but we can also use \emph{right-handed notation},
writing $\left(x\right)\!\overline{g}$ instead of $\overline{g}\!\left(x\right)$,
$\left(\left(x\right)\!\overline{h}\right)\!\overline{g}$ or $\left(x\right)\!\overline{h}\circ\overline{g}$
instead of $\overline{g}\!\left(\overline{h}\!\left(x\right)\right)$
or $\overline{g}\circ\overline{h}\!\left(x\right)$, and so on. A
convenient feature of right-handed notation is that the function first
applied to its argument comes first when reading from left to right.

To give expressions a smooth appearance, one may use \emph{overloaded}
notation, where function composition and function application are
denoted by the same symbol. For example, one may write $\overline{g}\circ\overline{h}\!\left(x\right)$
as $\overline{g}\circ\overline{h}\circ x$ and $\left(x\right)\!\overline{h}\circ\overline{g}$
as $x\circ\overline{h}\circ\overline{g}$ -- note that we do not need
parentheses due to the definition of function composition.

We may specifically refer to this type of notation as \emph{multiplicative}
overloaded notation. In certain contexts, it is more intuitive to
write $\overline{g}+\overline{h}+x$ or $x+\overline{h}+\overline{g}$,
using \emph{additive} overloaded notation instead of multiplicative. 

We can even overload all three operators involved, writing, for example,
\[
x+\overline{g}+\overline{h}=x+\overline{g+h}
\]
instead of $\left(x\right)\!\overline{g}\circ\overline{h}=\left(x\right)\!\overline{gh}$,
$x\circ\overline{g}\circ\overline{h}=x\circ\overline{gh}$ or $x+\overline{g}+\overline{h}=x+\overline{gh}$.
We call this \emph{totally overloaded notation}. 

It is useful to have different ways of writing $y/x$ (§ 5) in left-handed
and right-handed notation. We may write $y/x$ as $^{y}/_{x}$ when
using left-handed multiplicative notation and as $_{x}\setminus^{y}$
when using right-handed multiplicative notation; (\ref{eq:reg1})
then becomes $\overline{^{y}/_{x}}$$\,\circ\,x=y$ and $x\,\circ\,\overline{_{x}\setminus^{y}}=y$,
respectively.

The expression corresponding to $y/x$ in additive notation is $y-x$.
We may write this as $\left(y\shortleftarrow x\right)$ when using
left-handed notation and as $\left(x\shortrightarrow y\right)$ when
using right-handed notation. Relation (\ref{eq:reg1}) becomes $\overline{\left(y\shortleftarrow x\right)}+x=y$
and $x+\overline{\left(x\shortrightarrow y\right)}=y$, respectively;
both these expressions convey the idea that $x$ is moved to $y$
by some transformation on $X$.

\subsection{Covariant and contravariant actions of groups}

\paragraph*{14}

Closed actions of groups as defined above are \emph{covariant} in
the sense that $\overline{gh}\!\left(x\right)=\overline{g}\circ\overline{h}\!\left(x\right)$,
but it is also possible to define \emph{contravariant} actions such
that $\overline{gh}\!\left(x\right)=\overline{h}\circ\overline{g}\!\left(x\right)$.
In right-handed notation, we have $\left(x\right)\!\overline{gh}=\left(x\right)\!\overline{h}\circ\overline{g}$
for a covariant action, $\left(x\right)\!\overline{gh}=\left(x\right)\!\overline{g}\circ\overline{h}$
for a contravariant action. 

We note that covariant actions look most natural in left-handed notation,
whereas contravariant actions look most natural in right-handed notation.
As a consequence, covariant actions are traditionally called ``left
actions'', while contravariant actions are called ``right actions''.
Although the present terminology breaks with this tradition, it has
the virtue of clearly separating notation from concepts. (Furthermore,
I prefer the terminology covariant -- contravariant to alternatives
such as ``homomorphic'' -- ``antihomomorphic''.) 

It should be kept in mind, though, that covariant and contravariant
actions are exchangeable on the conceptual level. We can transform
a covariant (contravariant) action $\alpha$ of $G$ on $X$ to a
contravariant (covariant) action $\alpha^{\mathrm{op}}$ of $G$ on
$X$, defined by $\alpha^{\mathrm{op}}\!\left(gh\right):=\alpha\!\left(hg\right)$
for all $g,h\in G$.\newpage{}

\section{Sketch of an application:\protect \\
affine, (semi)preaffine and multiaffine spaces}

\subsection{Affine and preaffine spaces}

\paragraph*{15}

In this section, $V$ will denote a vector space over some field $K$,
with generic elements $\mathbf{u},\mathbf{v},\ldots$ and identity
element $\mathbf{0}$. An \emph{action }of $V$ on a set $X$ is an
action of the additive group of $V$ on $X$. We define an \emph{affine
space} as a non-empty set $X$ equipped with a \emph{regular monoidal}
action $\alpha$ of $V$ on $X$, or equivalently (see § 8) a regular
group action of $V$ on $X$ \cite{key-2}.

We shall use \emph{right-handed, totally overloaded additive notation}
in this section, so we write $\mathbf{u}\square\mathbf{v}$ as $\mathbf{u}+\mathbf{v}$,
$\overline{\mathbf{v}}\!\left(x\right)$ as $x+\overline{\mathbf{v}}$,
and $\overline{\mathbf{v}}\!\left(\overline{\mathbf{u}}\!\left(x\right)\right)$
as $x+\overline{\mathbf{u}}+\overline{\mathbf{v}}$. Closed actions
of groups are assumed to be \emph{contravariant} so that $\overline{\mathbf{u}+\mathbf{v}}\!\left(x\right)=\overline{\mathbf{v}}\!\left(\overline{\mathbf{u}}\!\left(x\right)\right)$,
written as $x+\overline{\mathbf{u}+\mathbf{v}}=x+\overline{\mathbf{u}}+\overline{\mathbf{v}}$.

By the definition of an affine space, $\alpha$ is a unital, closed
action of groups, so
\begin{gather}
x+\overline{\mathbf{0}}=x\qquad\forall x\in X,\label{eq:unitaladd}\\
\overline{\mathbf{u}}+\overline{\mathbf{v}}=\overline{\mathbf{u}+\mathbf{v}}\qquad\forall\mathbf{u},\mathbf{v}\in V.\label{eq:closedgrpadd}
\end{gather}
Furthermore, $\alpha$ is a regular action, so there is a function
\[
X\times X\rightarrow V,\qquad\left(x,y\right)\mapsto\left(x\shortrightarrow y\right)
\]
such that
\begin{gather}
x+\overline{\left(x\shortrightarrow y\right)}=y\qquad\forall x,y\in X,\label{eq:reg1add}\\
\left(x\shortrightarrow\left(x+\overline{\mathbf{v}}\right)\right)=\mathbf{v}\qquad\forall x\in X,\forall\mathbf{v}\in V,\label{eq:reg2add}
\end{gather}
corresponding to (\ref{eq:reg1}) and (\ref{eq:reg2}), respectively.
By (\ref{eq:reg1add}), $\mathbf{v}=\left(x\shortrightarrow y\right)$
is a solution of the equation $x+\overline{\mathbf{v}}=y$, and (\ref{eq:reg2add})
ensures that this solution is unique, because if $x+\overline{\mathbf{u}}=x+\overline{\mathbf{v}}$
then $\left(x\shortrightarrow\left(x+\overline{\mathbf{u}}\right)\right)=\left(x\shortrightarrow\left(x+\overline{\mathbf{v}}\right)\right)$,
so $\mathbf{u}=\mathbf{v}$.

The known properties of affine spaces can be derived from these assumptions.
For example, ``Chasles' law''
\begin{equation}
\left(x\shortrightarrow y\right)+\left(y\shortrightarrow z\right)=\left(x\shortrightarrow z\right)\qquad\forall x,y,z\in X,\label{eq:aChaslesV}
\end{equation}
can be derived by noting that (\ref{eq:closedgrpadd}) and (\ref{eq:reg1add})
imply 
\[
x+\overline{\left(x\shortrightarrow y\right)+\left(y\shortrightarrow z\right)}=x+\overline{\left(x\shortrightarrow y\right)}+\overline{\left(y\shortrightarrow z\right)}=z,
\]
and (\ref{eq:reg2add}) thus implies
\[
\left(x\shortrightarrow z\right)=\left(x\shortrightarrow\left(x+\overline{\left(x\shortrightarrow y\right)+\left(y\shortrightarrow z\right)}\right)\right)=\left(x\shortrightarrow y\right)+\left(y\shortrightarrow z\right).
\]

It is worth noting that
\begin{equation}
x+\overline{\mathbf{u}}+\overline{\mathbf{v}}=x+\overline{\mathbf{u}+\mathbf{v}}=x+\overline{\mathbf{v}+\mathbf{u}}=x+\overline{\mathbf{v}}+\overline{\mathbf{u}}\qquad\forall x\in X,\forall\mathbf{u},\mathbf{v}\in V.\label{eq:abelianT}
\end{equation}
Geometrically, this is the so-called parallelogram law of vector addition
in an affine space; see Figure 3.1.

\begin{figure}[h]
\begin{centering}
\includegraphics[scale=0.65]{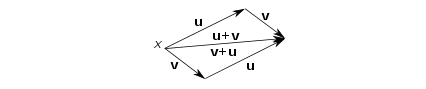}
\par\end{centering}

\caption{Composition of translations in an affine space.}
\end{figure}

\paragraph*{16}

We define a \emph{preaffine space} as a non-empty set $X$ equipped
with a \emph{regular premonoidal} action $\alpha$ of $V$ on $X$.
Identities (\ref{eq:unitaladd}), (\ref{eq:reg1add}) and (\ref{eq:reg2add})
are thus still required to hold, but instead of (\ref{eq:closedgrpadd})
we require that $\alpha$ is closed as an action of sets (§ 3). That
is, we assume that
\begin{equation}
\overline{\left(x\shortrightarrow y\right)}+\overline{\left(y\shortrightarrow z\right)}=\overline{\left(x\shortrightarrow z\right)}\qquad\forall x,y,z\in X.\label{eq:closedset}
\end{equation}
In other words, the transformation composed of the translation that
sends $x$ to $y$ followed by the translation that sends $y$ to
$z$ is itself a translation, necessarily the one that sends $x$
to $z$.

It is clear that an affine space is a preaffine space, since (\ref{eq:closedset})
follows from (\ref{eq:closedgrpadd}) and Chasles' law (\ref{eq:aChaslesV}).
A \emph{strictly preaffine} space is a preaffine space that is not
an affine space.\smallskip{}

\paragraph*{17}

In a strictly preaffine space, we do not have $\overline{\mathbf{u}}+\overline{\mathbf{v}}=\overline{\mathbf{u}+\mathbf{v}}$
for all $\mathbf{u},\mathbf{v}\in V$, and thus the parallelogram
law (\ref{eq:abelianT}) does not hold, and translations do not always
commute; see Figure 3.2.

\begin{figure}[h]
\begin{centering}
\includegraphics[scale=0.65]{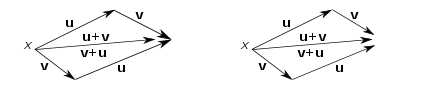}
\par\end{centering}

\caption{Commutative and non-commutative composition of translations in a non-affine
space.}
\end{figure}

On the other hand, recall from § 9 that the translations $\overline{\mathbf{v}}$
form a group, albeit not necessarily an abelian group. The unique
identity element in this group is $\overline{\mathbf{0}}$, and we
denote the unique inverse of $\overline{\mathbf{v}}$ by $-\overline{\mathbf{v}}$,
writing $\overline{\mathbf{u}}+\left(-\overline{\mathbf{v}}\right)$
as $\overline{\mathbf{u}}-\overline{\mathbf{v}}$. We have
\[
x+\overline{\left(x\shortrightarrow x\right)}=x=x+\overline{\mathbf{0}}\qquad\forall x\in X,
\]
and by § $5\dagger$ this implies $\left(x\shortrightarrow x\right)=\mathbf{0}$
and $\overline{\left(x\shortrightarrow x\right)}=\overline{\mathbf{0}}$
for every $x\in X$. Hence, $\overline{\left(x\shortrightarrow y\right)}+\overline{\left(y\shortrightarrow x\right)}=\overline{\mathbf{0}}$
by (\ref{eq:closedset}), so $\overline{\left(x\shortrightarrow y\right)}=\overline{\left(x\shortrightarrow y\right)}+\overline{\left(y\shortrightarrow x\right)}-\overline{\left(y\shortrightarrow x\right)}=-\overline{\left(y\shortrightarrow x\right)}$.

\subsection{Action fields, multiaffine spaces and associated semipreaffine spaces}

\paragraph*{18}

Recall that $x+\overline{\mathsf{g}}$ is another way of writing $\alpha\!\left(\mathsf{g}\right)\!\left(x\right)$.
Let us return to the latter notation temporarily while introducing
a generalization of the concept of action. An \emph{action field}
of $\mathsf{G}$ on $X$ is a function $\Phi$ that to every $p\in X$
assigns an action 
\[
\Phi\!\left(p\right):\mathsf{G}\rightarrow\mathscr{F}_{\!X}
\]
of $\mathsf{G}$ on $X$ in such a way that $\Phi\!\left(p\right)\!\left(\mathsf{G}\right)=\Phi\!\left(q\right)\!\left(\mathsf{G}\right)=:\overline{\mathsf{G}}$
for any $p,q\in X.$ When dealing with an action field, we thus have
expressions of the form $\Phi\!\left(p\right)\!\left(\mathsf{g}\right)\!\left(x\right)$
instead of $\alpha\!\left(\mathsf{g}\right)\!\left(x\right)$, or
$\Phi\!\left(p\right)\!\left(\mathbf{v}\right)\!\left(x\right)$ instead
of $\alpha\!\left(\mathsf{\mathbf{v}}\right)\!\left(x\right)$ if
$\mathsf{G}$ is the additive group of a vector space $V$. In the
latter case we of course write $\overline{V}$ instead of $\overline{\mathsf{G}}$.

In right-handed additive notation, we write $\Phi\!\left(p\right)\!\left(\mathbf{v}\right)\!\left(x\right)$
as
\begin{equation}
x+\overline{\mathbf{v}}^{p}.\label{eq:multi-pre1}
\end{equation}
 Analogously, we denote the unique solution $\mathbf{v}$ of the equation
$x+\overline{\mathbf{v}}^{p}=y$ by
\begin{equation}
\left(x\shortrightarrow y\right)_{p}.\label{eq:multipre2}
\end{equation}
For any $p\in X$, $\overline{\left(x\shortrightarrow y\right)_{p}}^{p}$
is thus the unique translation that sends $x$ to $y$. 

\smallskip{}

\paragraph*{19}

A \emph{regular monoidal action field} is an action field $\Phi$
of $V$ on $X$ such that $\Phi\!\left(p\right)$ is a regular monoidal
action of $V$ on $X$ for each $p\in X$. We call a non-empty set
$X$ equipped with a regular monoidal action field $\Phi$ a \emph{pointwise
affine} or \emph{multiaffine }space, since $X$ is an affine space
with respect to $\Phi\!\left(p\right)$ for each $p\in X$. (A \emph{strictly
multiaffine} space is a multiaffine space with a non-constant action
field.) Specifically, a multiaffine space is characterized by the
following identities:
\begin{gather}
x+\overline{\mathbf{0}}^{p}=x\qquad\forall p,x\in X,\label{eq:unitaladd-m}\\
x+\overline{\mathbf{u}}^{p}+\overline{\mathbf{v}}^{p}=x+\overline{\mathbf{u}+\mathbf{v}}^{p}\qquad\forall p,x\in X,\forall\mathbf{u},\mathbf{v}\in V,\label{eq:closedgrpadd-m}\\
x+\overline{\left(x\shortrightarrow y\right)_{p}}^{p}=y\qquad\forall p,x,y\in X,\label{eq:reg1add-m}\\
\left(x\shortrightarrow\left(x+\overline{\mathbf{v}}^{p}\right)\right)_{p}=\mathbf{v}\qquad\forall p,x\in X,\forall\mathbf{v}\in V.\label{eq:reg2add-m}
\end{gather}

A \emph{multipreaffine} space can be defined similarly as a non-empty
set $X$ equipped with a regular premonoidal action of $V$ on $X$
for each $p\in X$, but such spaces will not be considered here.\smallskip{}

\paragraph*{20}

If $X$ is a multiaffine space then each $\Phi\!\left(p\right)$ is
a free action, so by § $5\dagger$ there are for any $p\in X$ two
bijective functions,
\begin{gather*}
\phi\!\left(p\right):V\rightarrow\overline{V},\qquad\mathbf{v}\mapsto\phi\!\left(p\right)\!\left(\mathbf{v}\right)=\overline{\mathbf{v}}^{p},\\
\phi\!\left(p\right)^{-1}:\overline{V}\rightarrow V,\qquad\overline{\mathbf{v}}^{q}\mapsto\phi\!\left(p\right)^{-1}\!\left(\overline{\mathbf{v}}^{q}\right)=:\left(\overline{\mathbf{v}}^{q}\right)_{p}
\end{gather*}
such that $\phi\!\left(p\right)^{-1\!}\circ\phi\!\left(p\right)=\mathrm{Id}_{V}$
and $\phi\!\left(p\right)\circ\phi\!\left(p\right)^{-1\!}=\mathrm{Id}_{\overline{V}}$
for each $p\in X$. We also have $\phi\!\left(q\right)\circ\phi\!\left(q\right)^{-1\!}\circ\phi\!\left(p\right)=\phi\!\left(p\right)$
for any $p,q\in X$, and we can write
\begin{gather}
\left(\overline{\mathbf{v}}^{p}\right)_{p}=\mathbf{v}\qquad\forall p\in X,\forall\mathbf{v}\in V,\label{eq:inv1}\\
\overline{\left(\overline{\mathbf{v}}^{p}\right)_{q}}^{q}=\overline{\mathbf{v}}^{p}\qquad\forall p,q\in X,\forall\mathbf{v}\in V.\label{eq:inv1b}
\end{gather}
As $\overline{\left(x\shortrightarrow y\right)_{p}}^{p}=\overline{\left(x\shortrightarrow y\right)_{q}}^{q}$
for all $p,q\in X$, (\ref{eq:inv1}) gives 
\begin{equation}
\left(\overline{\left(x\shortrightarrow y\right)_{p}}^{p}\right)_{q}=\left(x\shortrightarrow y\right)_{q}\qquad\forall x,y,p,q\in X.\label{eq:inv2}
\end{equation}

\paragraph*{21}

Given a multiaffine space with action field $\Phi$ of $V$ on $X$,
we can define a single action $\alpha_{\Phi}$ of $V$ on $X$ by
setting $\alpha_{\Phi}\!\left(\mathbf{v}\right)\!\left(x\right)=\Phi\!\left(x\right)\!\left(\mathbf{v}\right)\!\left(x\right)$
for all $x\in X$ and $\mathbf{v}\in V$, or in the right-handed additive
notation used here, 
\begin{equation}
x+\overline{\mathbf{v}}:=x+\overline{\mathbf{v}}^{x}\qquad\forall x\in X,\forall\mathbf{v}\in V.\label{eq:redPreadd1}
\end{equation}
For consistency, we define 
\begin{equation}
\left(x\shortrightarrow y\right):=\left(x\shortrightarrow y\right)_{x}\qquad\forall x,y\in X.\label{eq:redPreadd2}
\end{equation}

Substituting $x$ for $p$ in (\ref{eq:unitaladd-m}), (\ref{eq:reg1add-m})
and (\ref{eq:reg2add-m}) and using (\ref{eq:redPreadd1}) and (\ref{eq:redPreadd2}),
we obtain (\ref{eq:unitaladd}), (\ref{eq:reg1add}) and (\ref{eq:reg2add}),
respectively. On the other hand, to obtain
\[
p+\overline{\left(x\shortrightarrow y\right)}+\overline{\left(y\shortrightarrow z\right)}=p+\overline{\left(x\shortrightarrow z\right)}\qquad\forall p,x,y,x\in X
\]
 by using (\ref{eq:redPreadd1}) and (\ref{eq:redPreadd2}) we need
to start from
\[
p+\overline{\left(x\shortrightarrow y\right)_{x}}^{p}+\overline{\left(y\shortrightarrow z\right)_{y}}^{p+\overline{\left(x\shortrightarrow y\right)_{x}}^{p}}=p+\overline{\left(x\shortrightarrow z\right)_{x}}^{p}\qquad\forall p,x,y,x\in X,
\]
but in general the latter relation holds only if $p=x$ (by \ref{eq:reg1add-m}).
This means that we cannot derive (\ref{eq:closedset}) from (\ref{eq:unitaladd-m})--(\ref{eq:reg2add-m})
by means of (\ref{eq:redPreadd1}) and (\ref{eq:redPreadd2}), and
\emph{a fortiori} (\ref{eq:closedgrpadd}) does not hold. We say that
$\alpha_{\Phi}$ defines a \emph{semipreaffine space}.

\pagebreak{}

\subsection{Deformation of affine spaces}

\paragraph*{22}

It is helpful to think of a non-affine space as an affine space which
has suffered some deformation. What separates an affine space from
a strictly preaffine space is whether or not the identity $\overline{\mathbf{u}}+\overline{\mathbf{v}}=\overline{\mathbf{u}+\mathbf{v}}$
holds, and in this section we will consider measures of ``non-affineness''
based on ``discrepancies'' such as that between $\overline{\mathbf{u}}+\overline{\mathbf{v}}$
and $\overline{\mathbf{u}+\mathbf{v}}$. We shall relate such discrepancies
to points in the space considered.

Let us first look at measures of the deformation at $x$ of an affine
space $X$, regarded as a non-strictly preaffine space, into a strictly
preaffine space. It is natural to start with the function $\mathfrak{T}_{0}:X\times V\times V\rightarrow V$,
defined by the equation
\begin{equation}
x+\overline{\mathbf{u}}+\overline{\mathbf{v}}+\overline{\mathfrak{T}_{0}\!\left(x,\mathbf{u},\mathbf{v}\right)}=x+\overline{\mathbf{u}+\mathbf{v}}\qquad\forall x\in X,\forall\mathbf{u},\mathbf{v}\in V\label{eq:simpledef}
\end{equation}
with the unique solution 
\begin{equation}
\mathfrak{T}_{0}\!\left(x,\mathbf{u},\mathbf{v}\right)=\left(\left(x+\overline{\mathbf{u}}+\overline{\mathbf{v}}\right)\shortrightarrow\left(x+\overline{\mathbf{u}+\mathbf{v}}\right)\right).\label{eq:t0}
\end{equation}
$\mathfrak{T}_{0}$ is a measure of ``general non-affineness'';
by definition $X$ is an affine space if and only if $\mathfrak{T}_{0}\!\left(x,\mathbf{u},\mathbf{v}\right)=\mathbf{0}$
for all $x\in X$ and $\mathbf{u},\mathbf{v}\in V$. It is obvious
that $\mathfrak{T}_{0}\!\left(x,\mathbf{0},\mathbf{v}\right)=\mathfrak{T}_{0}\!\left(x,\mathbf{u},\mathbf{0}\right)=\mathbf{0}$
for all $x\in X$ and all $\mathbf{u},\mathbf{v}\in V$. \smallskip{}

\paragraph*{23}

We shall now consider a special kind of deformation of affine spaces.
Modifying the definition in the previous paragraph, one can construct
a translation-skewsymmetric function $\mathfrak{T}_{1}:X\times V\times V\rightarrow V$
such that $\overline{\mathfrak{T}_{1}\!\left(x,\mathbf{u},\mathbf{v}\right)}=-\overline{\mathfrak{T}_{1}\!\left(x,\mathbf{v},\mathbf{u}\right)}$
by requiring that
\begin{equation}
x+\overline{\mathbf{u}}+\overline{\mathbf{v}}+\overline{\mathfrak{T}_{1}\!\left(x,\mathbf{u},\mathbf{v}\right)}=x+\overline{\mathbf{v}}+\overline{\mathbf{u}}\qquad\forall x\in X,\forall\mathbf{u},\mathbf{v}\in V,\label{eq:torseq-pre}
\end{equation}
 so that 
\begin{equation}
\mathfrak{T}_{1}\!\left(x,\mathbf{u},\mathbf{v}\right)=\left(\left(x+\overline{\mathbf{u}}+\overline{\mathbf{v}}\right)\shortrightarrow\left(x+\overline{\mathbf{v}}+\overline{\mathbf{u}}\right)\right).\label{eq:torsdef-pre}
\end{equation}
As $\left(x\rightarrow x\right)=\mathbf{0}$, $\mathfrak{T}_{1}\!\left(x,\mathbf{0},\mathbf{v}\right)=\mathfrak{\mathfrak{T}}_{1}\!\left(x,\mathbf{u},\mathbf{0}\right)=\mathbf{0}$
for all $x\in X$ and all $\mathbf{u},\mathbf{v}\in V$. Also, by
(\ref{eq:torsdef-pre}) and (\ref{eq:closedset}), $\overline{\mathfrak{T}_{1}\!\left(x,\mathbf{u},\mathbf{v}\right)}+\overline{\mathfrak{T}_{1}\!\left(x,\mathbf{v},\mathbf{u}\right)}=\overline{\left(\left(x+\overline{\mathbf{u}}+\overline{\mathbf{v}}\right)\shortrightarrow\left(x+\overline{\mathbf{u}}+\overline{\mathbf{v}}\right)\right)}=\overline{\mathbf{0}}$,
so $\overline{\mathfrak{T}_{1}\!\left(x,\mathbf{u},\mathbf{v}\right)}=-\overline{\mathfrak{T}_{1}\!\left(x,\mathbf{v},\mathbf{u}\right)}$
as promised. Furthermore,
\begin{gather*}
\overline{\mathfrak{T}_{0}\!\left(x,\mathbf{u},\mathbf{v}\right)}-\overline{\mathfrak{T}_{0}\!\left(x,\mathbf{v},\mathbf{u}\right)}=\overline{\left(\left(x+\overline{\mathbf{u}}+\overline{\mathbf{v}}\right)\shortrightarrow\left(x+\overline{\mathbf{u}+\mathbf{v}}\right)\right)}\\
+\overline{\left(\left(x+\overline{\mathbf{v}+\mathbf{u}}\right)\shortrightarrow\left(x+\overline{\mathbf{v}}+\overline{\mathbf{u}}\right)\right)}=\overline{\left(\left(x+\overline{\mathbf{u}}+\overline{\mathbf{v}}\right)\shortrightarrow\left(x+\overline{\mathbf{v}}+\overline{\mathbf{u}}\right)\right)},
\end{gather*}
so $\overline{\mathfrak{T}_{1}\!\left(x,\mathbf{u},\mathbf{v}\right)}=\overline{\mathfrak{T}_{0}\!\left(x,\mathbf{u},\mathbf{v}\right)}-\overline{\mathfrak{T}_{0}\!\left(x,\mathbf{v},\mathbf{u}\right)}$.

This relation between $\mathfrak{T}_{1}$ and $\mathfrak{T}_{0}$
shows that if $\mathfrak{T}_{0}\!\left(x,\mathbf{u},\mathbf{v}\right)=\mathbf{0}$
for all $x\in X$ and all $\mathbf{u},\mathbf{v}\in V$ then $\mathfrak{T}_{1}\!\left(x,\mathbf{u},\mathbf{v}\right)=\mathbf{0}$
for all $x\in X$ and all $\mathbf{u},\mathbf{v}\in V$. It is not
true, however, that $\mathfrak{T}_{1}\!\left(x,\mathbf{u},\mathbf{v}\right)=\mathbf{0}$
implies $\mathfrak{T}_{0}\!\left(x,\mathbf{u},\mathbf{v}\right)=\mathbf{0}$.
We say, informally at this point, that $\mathfrak{T}_{1}$ measures
the ``torsion'' of $X$, and that if $\mathfrak{T}_{1}\!\left(x,\mathbf{u},\mathbf{v}\right)=\mathbf{0}$
for all $x\in X$ and all $\mathbf{u},\mathbf{v}\in V$ then $X$
is ``torsion-free''. Thus, an affine space is torsion-free, but
a torsion-free space is not necessarily affine, as illustrated in
Figure 3.2. \smallskip{}

\paragraph*{24}

Let now $X$ be a multiaffine space associated with a semipreaffine
space. Instead of (\ref{eq:torseq-pre}), we have for the multiaffine
space
\begin{equation}
x+\overline{\mathbf{u}}^{x}+\overline{\mathbf{v}}^{y}+\overline{\mathfrak{T}_{1}^{*}\!\left(x,\mathbf{u},\mathbf{v}\right)}^{t}=x+\overline{\mathbf{v}}^{x}+\overline{\mathbf{u}}^{z}\qquad\forall x\in X,\forall\mathbf{u},\mathbf{v}\in V,\label{eq:mult-tors1-eq}
\end{equation}
where $x+\overline{\mathbf{u}}^{x}=y$, $y+\overline{\mathbf{v}}^{y}=t$,
and $x+\overline{\mathbf{v}}^{x}=z$, so 
\begin{equation}
\mathfrak{T}_{1}^{*}\!\left(x,\mathbf{u},\mathbf{v}\right)=\left(\left(x+\overline{\mathbf{u}}^{x}+\overline{\mathbf{v}}^{y}\right)\shortrightarrow\left(x+\overline{\mathbf{v}}^{x}+\overline{\mathbf{u}}^{z}\right)\right)_{t}.\label{eq:mult-tors1-fn}
\end{equation}
Corresponding to properties of $\mathfrak{T}_{1}$, we have $\mathfrak{T}_{1}^{*}\!\left(x,\mathbf{0},\mathbf{v}\right)=\mathfrak{T}_{1}^{*}\!\left(x,\mathbf{u},\mathbf{0}\right)=\mathbf{0}$
and $\overline{\mathfrak{T}_{1}^{*}\!\left(x,\mathbf{u},\mathbf{v}\right)}^{t}=-\overline{\mathfrak{T}_{1}^{*}\!\left(x,\mathbf{v},\mathbf{u}\right)}^{t}$
for all $x\in X$ and all $\mathbf{u},\mathbf{v}\in V$.

If $X$ has a constant action field, meaning that $X$ can be regarded
as an affine space, then $\overline{\mathbf{v}}^{y}=\overline{\mathbf{v}}^{x}$
and $\overline{\mathbf{u}}^{z}=\overline{\mathbf{u}}^{x}$, so $x+\overline{\mathbf{u}}^{x}+\overline{\mathbf{v}}^{y}=x+\overline{\mathbf{u}+\mathbf{v}}^{x}=x+\overline{\mathbf{v}+\mathbf{u}}^{x}=x+\overline{\mathbf{v}}^{x}+\overline{\mathbf{u}}^{z}$,
corresponding to (\ref{eq:abelianT}), so that $\mathfrak{T}_{1}^{*}\!\left(x,\mathbf{u},\mathbf{v}\right)=\mathbf{0}$
for all $x\in X$ and all $\mathbf{u},\mathbf{v}\in V$. This shows
again that an affine space is torsion-free. \smallskip{}

\paragraph*{25}

Let us continue with multiaffine spaces, defining $y$ and $z$ by
$x+\overline{\mathbf{u}}^{x}=y$ and $y+\overline{\mathbf{v}}^{y}=z$.
Corresponding to (\ref{eq:aChaslesV}), we have $\left(x\shortrightarrow y\right)_{p}+\left(y\shortrightarrow z\right)_{p}=\left(x\shortrightarrow z\right)_{p}$
for each $p\in X$, and this implies
\begin{equation}
x+\overline{\left(x\shortrightarrow y\right)_{x}}^{x}+\overline{\left(y\shortrightarrow z\right)_{y}}^{y}=x+\overline{\left(x\shortrightarrow y\right)_{x}+\left(y\shortrightarrow z\right)_{x}}^{x}\qquad\forall x,y,z\in X.\label{eq:curv000}
\end{equation}
We have $\left(x\shortrightarrow y\right)_{x}=\mathbf{u}$, $\left(y\shortrightarrow z\right)_{y}=\mathbf{v}$
and also $\left(\overline{\left(y\shortrightarrow z\right)_{y}}^{y}\right)_{x}=\left(y\shortrightarrow z\right)_{x}$
by (\ref{eq:inv2}), so we can rewrite (\ref{eq:curv000}) as
\begin{equation}
x+\overline{\mathbf{u}}^{x}+\overline{\mathbf{v}}^{y}=x+\overline{\mathbf{u}+\left(\overline{\mathbf{v}}^{y}\right)_{x}}^{x}\qquad\forall x\in X,\forall\mathbf{u},\mathbf{v}\in V.\label{eq:curv00}
\end{equation}

We can use (\ref{eq:curv00}) to obtain another measure of deformation
of affine spaces. Specifically, we inject $\overline{\mathbf{w}}^{x}$
into each side of the relation (\ref{eq:curv00}), obtaining
\begin{gather}
x+\overline{\mathbf{w}}^{x}+\overline{\mathbf{u}}^{r}+\overline{\mathbf{v}}^{s}=x+\overline{\mathbf{w}}^{x}+\overline{\mathbf{u}+\left(\overline{\mathbf{v}}^{s}\right)_{r}}^{r},\qquad\forall x\in X,\forall\mathbf{w},\mathbf{u},\mathbf{v}\in V,\label{eq:c0eq}
\end{gather}
where $x+\overline{\mathbf{w}}^{x}=r$ and $r+\overline{\mathbf{u}}^{r}=s$;
we also set $s+\overline{\mathbf{v}}^{s}=t$. Note that (\ref{eq:c0eq})
is a formal relation that does not necessarily hold; we can define
a function \linebreak{}
$\mathfrak{C}_{0}^{*}:X\times V\times V\times V\rightarrow V$ by
the equation 
\begin{equation}
x+\overline{\mathbf{w}}^{x}+\overline{\mathbf{u}}^{r}+\overline{\mathbf{v}}^{s}+\overline{\mathfrak{C}_{0}^{*}\left(x,\mathbf{w,}\mathbf{u},\mathbf{v}\right)}^{t}=x+\overline{\mathbf{w}}^{x}+\overline{\mathbf{u}+\left(\overline{\mathbf{v}}^{s}\right)_{r}}^{r}\label{eq:curveq}
\end{equation}
with the unique solution
\begin{equation}
\mathfrak{C}_{0}^{*}\!\left(x,\mathbf{w,}\mathbf{u},\mathbf{v}\right)=\left(\left(x+\overline{\mathbf{w}}^{x}+\overline{\mathbf{u}}^{r}+\overline{\mathbf{v}}^{s}\right)\shortrightarrow\left(x+\overline{\mathbf{w}}^{x}+\overline{\mathbf{u}+\left(\overline{\mathbf{v}}^{s}\right)_{r}}^{r}\right)\right)_{t}.\label{eq:curvdef}
\end{equation}

Clearly, $\mathfrak{C}_{0}^{*}\!\left(x,\mathbf{w},\mathbf{u},\mathbf{0}\right)\!=\!\mathbf{0}$
for all $x\in X$ and $\mathbf{w},\mathbf{u}\in V$, and by (\ref{eq:inv1b})
$\overline{\left(\overline{\mathbf{v}}^{s}\right)_{r}}^{r}\!=\!\overline{\mathbf{v}}^{s}$,
so $\mathfrak{C}_{0}^{*}\!\left(x,\mathbf{w,}\mathbf{0},\mathbf{v}\right)=\mathbf{0}$
for all $x\in X$ and $\mathbf{w},\mathbf{v}\in V$. Finally, if $\mathbf{w}=\mathbf{0}$
then (\ref{eq:c0eq}) reduces to (\ref{eq:curv00}), so $\mathfrak{C}_{0}^{*}\!\left(x,\mathbf{0},\mathbf{u},\mathbf{v}\right)=\mathbf{0}$
for all $x\in X$ and $\mathbf{u},\mathbf{v}\in V$.

If $X$ is a multiaffine space with a constant vector field so that
$X$ can be regarded as an affine space, we have $\left(\overline{\mathbf{v}}^{s}\right)_{r}=\left(\overline{\mathbf{v}}^{r}\right)_{r}=\mathbf{v}$
by (\ref{eq:inv1}). Adding to this the facts that $\overline{\mathbf{u}}^{r}=\overline{\mathbf{u}}^{x}$,
$\overline{\mathbf{v}}^{s}=\overline{\mathbf{v}}^{x}$, and $\overline{\mathbf{u}+\mathbf{v}}^{r}=\overline{\mathbf{u}+\mathbf{v}}^{x}=\overline{\mathbf{u}}^{x}+\overline{\mathbf{v}}^{x}$,
we conclude that (\ref{eq:c0eq}) holds. Thus, in an affine space
$\mathfrak{C}_{0}^{*}\!\left(x,\mathbf{w},\mathbf{u},\mathbf{v}\right)=\mathbf{0}$
for all $x\in X$ and $\mathbf{w},\mathbf{u},\mathbf{v}\in V$.

Without clear justification at the moment, we say that $\mathfrak{C}_{0}^{*}$
is a measure of ``curvature'', and $X$ is said to be ``flat''
if and only if $\mathfrak{C}_{0}^{*}\!\left(x,\mathbf{w},\mathbf{u},\mathbf{v}\right)=\mathbf{0}$
for all $x\in X$ and $\mathbf{w},\mathbf{u},\mathbf{v}\in V$. Thus,
an affine space is flat, but (\ref{eq:c0eq}) does not imply 
\[
x+\overline{\mathbf{u}}^{x}+\overline{\mathbf{v}}^{x+\overline{\mathbf{u}}^{x}}=x+\overline{\mathbf{u}+\mathbf{v}}^{x}\qquad\forall x\in X,\forall\mathbf{u},\mathbf{v}\in V,
\]
so a flat space is not necessarily affine. (Based on the displayed
relation, we could have defined a function $\mathfrak{T}_{0}^{*}$,
a measure of ``general non-affineness'' analogous to $\mathfrak{T}_{0}$.)\smallskip{}

\paragraph*{26}

One obtains a function $\mathfrak{C}_{1}^{*}:X\times V\times V\times V\rightarrow V$
corresponding to $\mathfrak{C}_{0}^{*}$ and translation-skewsymmetric
in $\mathbf{u}$ and $\mathbf{v}$ by setting 
\begin{equation}
\overline{\mathfrak{C}_{1}^{*}\!\left(x,\mathbf{w,}\mathbf{u},\mathbf{v}\right)}=\overline{\mathfrak{C}_{0}^{*}\!\left(x,\mathbf{w,}\mathbf{u},\mathbf{v}\right)}-\overline{\mathfrak{C}_{0}^{*}\!\left(x,\mathbf{w,}\mathbf{v},\mathbf{u}\right)}.\label{eq:c1def}
\end{equation}
Reasoning as in § 23, we obtain $\overline{\mathfrak{C}_{1}^{*}\!\left(x,\mathbf{w,}\mathbf{u},\mathbf{v}\right)}=-\overline{\mathfrak{C}_{1}^{*}\!\left(x,\mathbf{w,}\mathbf{v},\mathbf{u}\right)}$.
We also have $\mathfrak{C}_{1}^{*}\!\left(x,\mathbf{\mathbf{0},u},\mathbf{v}\right)\!=\mathfrak{\mathfrak{C}}_{1}^{\ast}\!\left(x,\mathbf{\mathbf{w},0},\mathbf{v}\right)\!=\mathfrak{\mathfrak{C}_{1}^{*}}\!\left(x,\mathbf{w},\mathbf{u},\mathbf{0}\right)=\mathbf{0}$
for all $x\in X$ and all $\mathbf{w},\mathbf{u},\mathbf{v}\in V$.\pagebreak{}

\subsection{Bound vectors and parallel transport}

\paragraph*{27}

The measures of deformation of an affine space derived in the previous
subsection can be interpreted geometrically in terms of bound vectors
and parallel transport.

A \emph{bound vector} is defined for present purposes simply as a
pair of points $\left(x,y\right),$ intuitively corresponding to a
directed line segment connecting $x$ to $y$. We say that $x$ is
the \emph{origin} and $y$ the \emph{tip} of $\left(x,y\right)$.

\emph{Parallel transport }of a bound vector translates its origin
and tip by means of the same vector. This is often described by speaking
about parallel transport of a bound vector \emph{along} another bound
vector with the same origin. 

In a preaffine space, parallel transport by $\mathbf{v}$, or along
$\left(x,x\!+\!\overline{\mathbf{v}}\right)$, maps $\left(x,x\!+\!\overline{\mathbf{u}}\right)$
to $\left(x+\overline{\mathbf{v}},x+\overline{\mathbf{u}}+\overline{\mathbf{v}}\right)$.
Alternatively, parallel transport along $\left(x,y\right)$ maps $\left(x,z\right)$
to $\left(x+\overline{\left(x\shortrightarrow y\right)},z+\overline{\left(x\shortrightarrow y\right)}\right)=\left(y,z+\overline{\left(x\shortrightarrow y\right)}\right)$.

In a multiaffine space, parallel transport along $\left(x,x+\overline{\mathbf{v}}^{x}\right)$
maps $\left(x,x+\overline{\mathbf{u}}^{x}\right)$ to $\left(x+\overline{\mathbf{v}}^{x},x+\overline{\mathbf{u}}^{x}+\overline{\mathbf{v}}^{y}\right)$,
where $y=x+\overline{\mathbf{u}}^{x}$. Alternatively, parallel transport
along $\left(x,y\right)$ maps $\left(x,z\right)$ to $\left(x+\overline{\left(x\shortrightarrow y\right)_{x}}^{x},z+\overline{\left(x\shortrightarrow y\right)_{x}}^{y}\right)=\left(y,z+\overline{\left(x\shortrightarrow y\right)_{x}}^{y}\right)$.
\smallskip{}

\paragraph*{28}

Let the points $x,y,z\in X$ be given. Map $\left(x,z\right)$ to
$\left(y,t\right)$ by parallel transport along $\left(x,y\right)$,
and map $\left(x,y\right)$ to $\left(z,t'\right)$ by parallel transport
along $\left(x,z\right)$. $\mathfrak{T}_{1}$, or equivalently $\mathfrak{T}_{1}^{*}$,
is a measure of the discrepancy between the tips $t=x+\overline{\mathbf{u}}+\mathbf{\overline{\mathbf{v}}}$
(or $t=x+\overline{\mathbf{u}}^{x}+\mathbf{\overline{\mathbf{v}}^{y}}$)
and $t'=x+\overline{\mathbf{v}}+\mathbf{\overline{\mathbf{u}}}$ (or
$t=x+\overline{\mathbf{v}}^{x}+\mathbf{\overline{\mathbf{u}}}^{z}$)
of the two parallel-transported bound vectors $\left(y,t\right)$
and $\left(z,t'\right)$. Figure 3.3 visualizes this parallel-transport
construction.

\begin{figure}[h]
\begin{centering}
\includegraphics[scale=0.65]{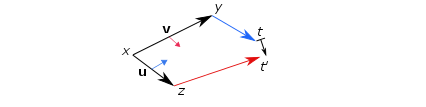}
\par\end{centering}

\caption{Geometrical interpretation of $\mathfrak{T}_{1}$ or $\mathfrak{T}_{1}^{*}$.}
\end{figure}

\paragraph*{29}

Let the points $x,y,z,r\in X$ be given. Map $\left(x,r\right)$ to
$\left(y,s\right)$ by parallel transport along $\left(x,y\right)$,
and then map $\left(y,s\right)$ to $\left(z,t\right)$ by parallel
transport along $\left(y,z\right)$. Also map $\left(x,r\right)$
to $\left(z,t'\right)$ by parallel transport along $\left(x,z\right)$.
$\mathfrak{C}_{0}^{*}$ is a measure of the discrepancy between the
tips $t=x+\overline{\mathbf{w}}^{x}+\overline{\mathbf{u}}^{r}+\overline{\mathbf{v}}^{s}$
and $t'=x+\overline{\mathbf{w}}^{x}+\overline{\mathbf{u}+\left(\overline{\mathbf{v}}^{s}\right)_{r}}^{r}$
of the parallel-transported bound vectors $\left(z,t\right)$ and
$\left(z,t'\right)$, where $z=x+\overline{\mathbf{u}}^{x}+\overline{\mathbf{v}}^{y}=x+\overline{\mathbf{u}+\left(\overline{\mathbf{v}}^{y}\right)_{x}}^{x}$.
Figure 3.4 shows this parallel-transport construction.

\begin{figure}[h]
\begin{centering}
\includegraphics[scale=0.65]{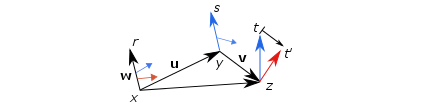}
\par\end{centering}

\caption{Geometrical interpretation of $\mathfrak{C}_{0}^{*}$.}
\end{figure}
\pagebreak{}

\subsection{Final remarks}

\paragraph*{30}

Readers may have suspected that what is presented here as an application
of a theory was historically a problem that helped to inspire the
construction of that theory. It is clear -- for example, from synthetic
differential geometry \cite{key-4} -- that the measures $\mathfrak{T}_{1}$
and $\mathfrak{T}_{1}^{*}$ are related to the torsion tensor in differential
geometry, while $\mathfrak{C}_{0}^{*}$ and $\mathfrak{C}_{1}^{*}$
are related to the curvature tensor. Note, though, that the present
notions of torsion and curvature are independent of the ``differential''
part of differential geometry -- it is possible to explain and understand
the general idea of parallel transport, torsion and curvature without
first introducing infinitesimal notions (of various orders). These
notions can be introduced at the next stage instead, generalizing
generalized affine spaces. 

The present approach to ``finitary differential geometry'' not only
opens a fast track to torsion and curvature but also suggests a path
from coordinate-free affine geometry defined in terms of group actions
to coordinate-free differential geometry. Recall that a \emph{pointwise
affine} space may yet be \emph{globally non-affine.} In other words,
$x+\overline{\mathbf{u}}^{p}+\overline{\mathbf{v}}^{p}=x+\overline{\mathbf{u}+\mathbf{v}}^{p}$
for all $p,x\in X$ and $\mathbf{u},\mathbf{v}\in V$ does not imply
\linebreak{}
$x+\overline{\mathbf{u}}^{x}+\overline{\mathbf{v}}^{\left(x+\overline{\mathbf{u}}^{x}\right)}=x+\overline{\mathbf{u}+\mathbf{v}}^{x}$
for all $x\in X$ and $\mathbf{u},\mathbf{v}\in V$ . In differential
geometry, on the other hand, we typically deal with a space that is
(1) \emph{pointwise} and (2) \emph{infinitesimally} affine but (3)
\emph{locally} and (4) \emph{globally} non-affine. That is, 
\begin{enumerate}
\item $x+\overline{\mathbf{u}}^{p}+\overline{\mathbf{v}}^{p}=x+\overline{\mathbf{u}+\mathbf{v}}^{p}$
for all $p,x\in X$ and $\mathbf{u},\mathbf{v}\in V$;
\item $x+\overline{\mathbf{u}}^{x}+\overline{\mathbf{v}}^{\left(x+\overline{\mathbf{u}}^{x}\right)}=x+\overline{\mathbf{u}+\mathbf{v}}^{x}$
for all $x\in X$ and all ``infinitesimal'' $\mathbf{u},\mathbf{v}\in V$
-- $X$ is ``infinitesimally affine'';
\item we do not have $x+\overline{\mathbf{u}}^{x}+\overline{\mathbf{v}}^{\left(x+\overline{\mathbf{u}}^{x}\right)}=x+\overline{\mathbf{u}+\mathbf{v}}^{x}$
for all $x\in X$ and all ``nearly infinitesimal'' $\mathbf{u},\mathbf{v}\in V$,
but $\mathfrak{T}_{0}^{*}$, given by
\[
\mathfrak{T}_{0}^{*}\!\left(x,\mathbf{u},\mathbf{v}\right)=\left(x+\overline{\mathbf{u}}^{x}+\overline{\mathbf{v}}^{\left(x+\overline{\mathbf{u}}^{x}\right)}\shortrightarrow x+\overline{\mathbf{u}+\mathbf{v}}^{x}\right)_{x+\overline{\mathbf{u}}^{x}+\overline{\mathbf{v}}^{\left(x+\overline{\mathbf{u}}^{x}\right)}},
\]
 is a bilinear function of its ``nearly infinitesimal'' vector arguments
-- $X$ is ``locally linearly connected'';
\item we do not have $x+\overline{\mathbf{u}}^{x}+\overline{\mathbf{v}}^{\left(x+\overline{\mathbf{u}}^{x}\right)}=x+\overline{\mathbf{u}+\mathbf{v}}^{x}$
for all $x\in X$ and all $\mathbf{u},\mathbf{v}\in V$.
\end{enumerate}
It remains to give formal definitions of notions such as ``infinitesimal''
and ``nearly infinitesimal''; this is only a heuristic description.
Note, though, that the values of $\mathfrak{T}_{0}$, $\mathfrak{T}_{0}^{*}$,
$\mathfrak{T}_{1}$, $\mathfrak{T}_{1}^{*}$, $\mathfrak{C}_{0}^{*}$,
and $\mathfrak{C}_{1}^{*}$ are all equal to $\mathbf{0}$ when one
of their vector arguments is $\mathbf{0}$, as required if these functions
are to be approximated by linear functions locally. In addition, $\mathfrak{D}^{*}$,
defined by
\[
x+\overline{\mathbf{v}}^{p}+\overline{\mathfrak{D}^{*}\!\left(x,p,\mathbf{d},\mathbf{v}\right)}^{p}=x+\overline{\mathbf{v}}^{p+\overline{\mathbf{d}}^{p}},
\]
satisfies $\mathfrak{D}^{*}\!\left(x,p,\mathbf{0},\mathbf{v}\right)=\mathbf{0}$
for all $x,p\in X$ and $\mathbf{v}\in V$, establishing a necessary
condition for $\mathfrak{D}^{*}$ to be approximated locally (for
$\mathbf{d}\approx\mathbf{0}$) by a linear function of $\mathbf{d}$,
and $x+\overline{\mathbf{v}}^{p'}$ to be approximated locally (for
$p'\approx p$) by $x+\overline{\mathbf{v}}^{p}+\overline{\mathfrak{D}^{*}\!\left(x,p,\left(p\shortrightarrow p'\right)_{p},\mathbf{v}\right)}^{p}$,
so that we have a ``smooth'' action field. There would seem to be
enough similarities between the finitary and the general cases to
make a generalization possible, yet also enough differences to make
that generalization nontrivial.\newpage{}

\appendix

\section{Malcev operations and related spaces}

\paragraph*{A1}

Let $\alpha$ be a regular premonoidal action of a vector space $V$
on a non-empty set $X$ (see § 16). We can define a ternary operation
\[
\kappa:X\times X\times X\rightarrow X,\qquad\left(x,y,z\right)\mapsto\left[x,y,z\right]:=x+\overline{\left(y\shortrightarrow z\right)}.
\]
 Let us express the postulates that define a \emph{preaffine }space
in terms of $\kappa$.

As $\alpha$ is a unital action, we have $x+\overline{0}=x$ for all
$x\in X$ (\ref{eq:unitaladd}), or equivalently $x+\overline{y\shortrightarrow y}=x$
for all $x,y\in X$. That is, we have the identity
\begin{equation}
\left[x,y,y\right]=x.\label{eq:xyy}
\end{equation}
Furthermore, $\alpha$ is a transitive action, so we have $x+\overline{x\shortrightarrow y}=y$
for all $x,y\in X$ (\ref{eq:reg1add}). Thus, we have the identity
\begin{equation}
\left[x,x,y\right]=y.\label{eq:xxy}
\end{equation}
 A ternary operation satisfying (\ref{eq:xyy}) and (\ref{eq:xxy})
is called a \emph{Malcev operation} \cite{key-3}. 

As $\alpha$ is a free action as well, we have $\left(x\shortrightarrow\left(x+\overline{\mathbf{v}}\right)\right)=\mathbf{v}$
(\ref{eq:reg2add}) for all $x\in X$, so $\left(x\shortrightarrow\left(x+\overline{\left(y\shortrightarrow z\right)}\right)\right)=\left(y\shortrightarrow z\right)$
for all $x,y,z\in X$, so $p+\overline{\left(x\shortrightarrow\left(x+\overline{\left(y\shortrightarrow z\right)}\right)\right)}=p+\overline{\left(y\shortrightarrow z\right)}$
for all $p,x,y,z\in X$. In terms of $\kappa$, this is the identity
\begin{equation}
\left[p,x,\left[x,y,z\right]\right]=\left[p,y,z\right].\label{eq:pxxyz}
\end{equation}

Finally, as $\alpha$ is a closed action of sets, we have $\overline{\left(x\shortrightarrow y\right)}+\overline{\left(y\shortrightarrow z\right)}=\overline{\left(x\shortrightarrow z\right)}$
(\ref{eq:closedset}) for all $x,y,x\in X$. Hence, $\left(p+\overline{\left(x\shortrightarrow y\right)}\right)+\overline{\left(y\shortrightarrow z\right)}=p+\overline{\left(x\shortrightarrow z\right)}$
for all $p,x,y,x\in X$. This translates to the identity
\begin{equation}
\left[\left[p,x,y\right],y,z\right]=\left[p,x,z\right].\label{eq:pxyyz}
\end{equation}

In a \emph{semipreaffine }space (\ref{eq:xyy}), (\ref{eq:xxy}),
and (\ref{eq:pxxyz}) hold; a \emph{strictly} semipreaffine space
is one where (\ref{eq:pxyyz}) does not hold.

$\kappa$ is said to be \emph{commutative} when $\left[x,y,z\right]=\left[z,y,x\right]$
for all $x,y,z\in X$, and \emph{associative} when $\left[\left[x,y,z\right],r,t\right]=\left[x,y,\left[z,r,t\right]\right]$
for all $x,y,z,r,t\in X$. A space with an associative Malcev operation
has many names, including \emph{heap} and \emph{groud}. A Malcev operation
satisfying either (\ref{eq:pxxyz}) or (\ref{eq:pxyyz}) is sometimes
said to be \emph{semiassociative} \cite{key-3} Note that if $\kappa$
is associative then (\ref{eq:xyy}) implies (\ref{eq:pxxyz}) and
(\ref{eq:xxy}) implies (\ref{eq:pxyyz}). Conversely, if (\ref{eq:pxxyz})
and (\ref{eq:pxyyz}) hold then
\[
\left[\left[x,y,z\right],r,t\right]=\left[\left[x,y,z\right],z,\left[z,r,t\right]\right]=\left[x,y,\left[z,r,t\right]\right],
\]
 so a Malcev operation is associative if and only if (\ref{eq:pxxyz})
and (\ref{eq:pxyyz}) hold.

It is easy to verify that a preaffine space is torsion-free (as described
in § 23 and §~28) if and only if $\left[y,x,z\right]=\left[z,x,y\right]$
for all $x,y,z\in X$, and that a semipreaffine space is flat (as
described in § 25 and § 29) if and only if $\left[\left[p,x,y\right],y,z\right]=\left[p,x,z\right]$
for all $p,x,y,z\in X$, that is, if (\ref{eq:pxyyz}) holds so that
$X$ is in fact a preaffine space. Thus, a commutative $\kappa$ gives
a torsion-free space, while an associative $\kappa$ gives a curvature-free
space. In other words, we can have non-zero curvature only in a strictly
semipreaffine space (or a corresponding strictly multiaffine space).

\smallskip{}
\texttt{ }

\paragraph*{A2}

Kock \cite{key-4} defines an \emph{affine connection} on $X$ as
a ternary operation $\lambda$ on $X$ such that $\lambda\left(x,x,y\right)=y$
and $\lambda\left(x,y,x\right)=y$ for all $x,y\in X$. (We disregard
here the ``infinitesimal'' part of Kock's argument.) Writing $\lambda\left(x,y,z\right)$
as $\left[y,x,z\right]$, these identities are rendered as $\left[x,x,y\right]=y$
(\ref{eq:xxy}) and $\left[y,x,x\right]=y$ (\ref{eq:xyy}), so $\lambda$
is in effect a Malcev operation. Kock then proves (using further assumptions)
that $\lambda\left(y,x,\lambda\left(x,y,z\right)\right)=z$ and $\lambda\left(z,\lambda\left(x,y,z\right),x\right)=y$,
or equivalently 
\begin{equation}
\left[x,y,\left[y,x,z\right]\right]=z\label{eq:k3}
\end{equation}
 and 
\begin{equation}
\left[\left[y,x,z\right],z,x\right]=y.\label{eq:k4}
\end{equation}
 It is clear that (\ref{eq:k3}) is implied by (A.2) and (A.3) while
(\ref{eq:k4}) is implied by (A.1) and (A.4). On the other hand, a
``Kock space'' is not necessarily a preaffine space. As we shall
see, it is sometimes possible to replace (A.3) by (\ref{eq:k3}) and
(A.4) by (\ref{eq:k4}) without losing desired properties of the spaces
considered.

\medskip{}

\paragraph*{A3}

In Section 3, we defined $x+\overline{\left(y\shortrightarrow z\right)}$
in terms of an action of a vector space on $X$. One can think of
$\left[x,y,z\right]$ as just a notational shorthand for $x+\overline{\left(y\shortrightarrow z\right)}$.
On the other hand, one can regard $\kappa$, required to satisfy (\ref{eq:xyy}),
(\ref{eq:xxy}) and one or more of (A.3)--(A.6), as a Malcev operation
that makes it possible to define preaffine and other spaces \emph{intrinsically},
without reference to an action of an external vector space. Thus,
a preaffine space is given by an abstract ternary operator satisfying
(A.1)--(A.4), a semipreaffine space is given by an abstract ternary
operator satisfying (A.1)--(A.3), and so on.

In terms of $\kappa$, a \emph{translation} on $X$ is a function
\begin{equation}
\left[-,a,b\right]:X\rightarrow X,\qquad x\mapsto\left[x,a,b\right]\qquad a,b\in X.\label{eq:deftrans}
\end{equation}
The set of all translations of this form is denoted $\mathscr{T}_{X}$.
As $\left[a,a,b\right]=b$ for all $a,b\in X$ by (\ref{eq:xxy}),
we have $\left[-,a,b\right]\!\left(a\right)=b$. For any $a\in X$,
the translation $\left[-,a,a\right]$ is the identity map $\epsilon_{X}$
on $X$ because $\left[x,a,a\right]=x$ for any $x,a\in X$ by (\ref{eq:xyy}).
Conversely, if a transformation of the form $\left[-,a,b\right]$
is $\epsilon_{X}$ so that $\left[-,a,b\right]\left(a\right)=a$ then
$b=\left[a,a,b\right]=a$ by (\ref{eq:xxy}).

Note that $\left[-,a,b\right]=\left[-,a',\left[a',a,b\right]\right]$
by (\ref{eq:pxxyz}). Conversely if $\left[-,a',b'\right]=\left[-,a,b\right]$
so that $\left[a',a',b'\right]=\left[a',a,b\right]$ then $b'=\left[a',a,b\right]$
by (\ref{eq:xxy}). Thus,
\[
\left[-,a',b'\right]=\left[-,a,b\right]\quad\Longleftrightarrow\quad b'=\left[a',a,b\right].
\]

Also note that if $\left[-,a,b\right]\left(p\right)=p$ for some $p\in X$
so that $\left[p,a,b\right]=p$ then
\[
\left[-,a,b\right]\left(x\right)=\left[x,a,b\right]=\left[x,p,\left[p,a,b\right]\right]=\left(x,p,p\right)=x
\]
for any $x\in X$ by (\ref{eq:pxxyz}) and (\ref{eq:xyy}), meaning
that if $\left[-,a,b\right]$ has a fixed point then $\left[-,a,b\right]=\epsilon_{X}$.
Thus translations defined in terms of $\kappa$ has this property
in common with translations defined in terms of a group actions (see
§ 9). In addition, if (\ref{eq:k4}) holds then $\left[-,a,b\right]\circ\left[-,b,a\right]=\epsilon_{X}$,
because 
\[
\left[-,a,b\right]\left(\left[-,b,a\right]\left(x\right)\right)=\left[\left[x,b,a\right],a,b\right]=x,
\]
so a translation as defined in (\ref{eq:deftrans}) is a bijective
function.

By definition, $\left[-,c,d\right]\circ\left[-,a,b\right]\left(x\right)=\left[-,c,d\right]\left(\left[-,a,b\right]\left(x\right)\right)=\left[\left[x,a,b\right],c,d\right]$.
Thus, the composite, called the \emph{sum}, of $\left[-,a,b\right]$
and $\left[-,c,d\right]$ is the transformation
\begin{gather}
\left[-,a,b\right]+\left[-,c,d\right]:X\rightarrow X\label{eq:defaddtla}\\
\left(x\mapsto\left[\left[x,a,b\right],c,d\right]\right)=:\left[\left[-,a,b\right],c,d\right].\nonumber 
\end{gather}
This definition is legitimate when (\ref{eq:xxy}) and (\ref{eq:pxxyz})
hold, because then $\left[-,a',b'\right]=\left[-,a,b\right]$ and
$\left[-,c',d'\right]=\left[-,c,d\right]$ implies 
\begin{gather*}
\left[\left[x,a',b'\right],c',d'\right]=\left[\left[x,a',\left[a',a,b\right]\right],c',\left[c',c,d\right]\right]=\left[\left[x,a,b\right],c,d\right].
\end{gather*}

To extend the definition of addition of translations, we set
\begin{gather*}
\left[-,a_{1},b_{1}\right]+\left(\left[-,a_{2},b_{2}\right]+\left[-,a_{3},b_{3}\right]\right):=\left(\left[-,a_{2},b_{2}\right]+\left[-,a_{3},b_{3}\right]\right)\circ\left[-,a_{1},b_{1}\right];\\
\left(\left[-,a_{1},b_{1}\right]+\left[-,a_{2},b_{2}\right]\right)+\left[-,a_{3},b_{3}\right]:=\left[-,a_{3},b_{3}\right]\circ\left(\left[-,a_{1},b_{1}\right]+\left[-,a_{2},b_{2}\right]\right).
\end{gather*}
Then
\begin{gather*}
\left[-,a_{1},b_{1}\right]+\left(\left[-,a_{2},b_{2}\right]+\left[-,a_{3},b_{3}\right]\right)=\left(\left[-,a_{3},b_{3}\right]\circ\left[-,a_{2},b_{2}\right]\right)\circ\left[-,a_{1},b_{1}\right]\\
=\left[-,a_{3},b_{3}\right]\circ\left(\left[-,a_{2},b_{2}\right]\circ\left[-,a_{1},b_{1}\right]\right)=\left(\left[-,a_{1},b_{1}\right]+\left[-,a_{2},b_{2}\right]\right)+\left[-,a_{3},b_{3}\right].
\end{gather*}
We may thus omit parentheses from sums involving translations, since
all ways of completing an expression of the form 
\[
\left[-,a_{1},b_{1}\right]+\left[-,a_{2},b_{2}\right]+\ldots+\left[-,a_{n},b_{n}\right]
\]
with parentheses yield the same transformation. We call a transformation
of this kind an \emph{iteration of translations}, and denote the set
of all iterations of translations on $X$ by $\mathscr{T}_{X}^{*}.$
Translations are thus a special case of iterations of translations,
and addition of translations is a special case of addition of iterations
of translations. Note that
\begin{gather*}
\left[-,a_{1},b_{1}\right]+\left[-,a_{2},b_{2}\right]+\ldots+\left[-,a_{n},b_{n}\right]\\
=\left[-,a_{n},b_{n}\right]\circ\ldots\circ\left[-,a_{2},b_{2}\right]\circ\left[-,a_{1},b_{1}\right]\\
=\left[\left[\left[\left[-,a_{1},b_{1}\right],a_{2},b_{2}\right],\ldots\right],a_{n},b_{n}\right],
\end{gather*}
so the sum of iterations of translations again does not depend on
the choice of pairs of points representing the translations involved.

It is clear that $\left(\mathscr{T}_{X}^{*},+\right)$, where + denotes
addition of iterations of translations, is an associative magma.

By (\ref{eq:xyy}), $\left[-,e,e\right]\!\left(x\right)=\left[x,e,e\right]=x$
for any $x,e\in X$, so $\left[-,e,e\right]=\epsilon_{X}$ for all
$e\in X$. By (\ref{eq:xyy}), we also have 
\begin{gather*}
\left(\left[-,e,e\right]+\left[-,c,d\right]\right)\!\left(x\right)=\left[\left[-,e,e\right]\!,c,d\right]\!\left(x\right)=\left[\left[x,e,e\right]\!,c,d\right]=\left[x,c,d\right]=\left[-,c,d\right]\!\left(x\right);\\
\left(\left[-,a,b\right]+\left[-,e,e\right]\right)\!\left(x\right)=\left[\left[-,a,b\right]\!,e,e\right]\!\left(x\right)=\left[\left[x,a,b\right]\!,e,e\right]=\left[x,a,b\right]=\left[-,a,b\right]\!\left(x\right).
\end{gather*}
Thus, for any $e\in X$ we have 
\begin{gather*}
\left[-,a,b\right]+\left[-,e,e\right]=\left[-,e,e\right]+\left[-,a,b\right]=\left[-,a,b\right];\\
\left[-,a_{1},b_{1}\right]+\ldots+\left[-,a_{n},b_{n}\right]+\left[-,e,e\right]=\left[-,e,e\right]+\left[-,a_{1},b_{1}\right]+\ldots+\left[-,a_{n},b_{n}\right]\\
=\left[-,a_{1},b_{1}\right]+\ldots+\left[-,a_{n},b_{n}\right]
\end{gather*}
 for any $\left[-,a,b\right]\in\mathscr{T}_{X}$ and $\left[-,a_{1},b_{1}\right]+\ldots+\left[-,a_{n},b_{n}\right]\in\mathscr{T}_{X}^{*}$.
Hence, $\left(\mathscr{T}_{X}^{*},+\right)$ is a monoid, since it
is an associative magma.

Furthermore, if (\ref{eq:k4}) holds then 
\[
\left(\left[-,a,b\right]+\left[-,b,a\right]\right)\left(x\right)=\left[\left[-,a,b\right],b,a\right]\left(x\right)=\left[\left[x,a,b\right],b,a\right]=x,
\]
so 
\[
\left[-,a,b\right]+\left[-,b,a\right]=\epsilon_{X}=\left[-,e,e\right].
\]
 Thus, every iteration of translations $\left[-,a_{1},b_{1}\right]+\ldots+\left[-,a_{n},b_{n}\right]$
has an inverse iteration of translations $\left[-,b_{n},a_{n}\right]+\ldots+\left[-,b_{1},a_{1}\right]$,
so $\left(\mathscr{T}_{X}^{*},+\right)$ is a group since it is a
monoid. In particular, if $\kappa$ is associative then $\left(\mathscr{T}_{X}^{*},+\right)$
is a group. 

If $\kappa$ is associative, we also have 
\begin{equation}
\left[-,a,b\right]+\left[-,c,d\right]=\left(x\mapsto\left[\left[x,a,b\right],c,d\right]\right)=\left(x\mapsto\left[x,a,\left[b,c,d\right]\right]\right).\label{eq:defaddtlb}
\end{equation}
Hence, $\left[-,a,b\right]+\left[-,c,d\right]$ is a translation as
well, and $\left(\mathscr{T}_{X},+\right)$ is a group because $\left(\mathscr{T}_{X}^{*},+\right)$
is a group (compare § 9).

\medskip{}

\paragraph*{A4}

Fix some $e\in X$. We call a translation that can be written in the
form $\left[-,e,b\right]$ a \emph{pointed translation} relative to
$e$. Any translation can be written as a pointed translation in at
most one way since 
\[
\left[-,e,b\right]=\left[-,e,b'\right]\quad\Longleftrightarrow\quad b=b'
\]
by (\ref{eq:xxy}). Hence, we can define the \emph{pointed sum} of
the pointed translations $\left[-,e,b\right]$ and $\left[-,e,d\right]$
as the pointed translation
\begin{gather}
\left[-,e,b\right]+_{e}\left[-,e,d\right]:X\rightarrow X,\label{eq:defaddtl2}\\
\left(x\mapsto\left[x,e,\left[b,e,d\right]\right]\right)=:\left[-,e,\left[b,e,d\right]\right].\nonumber 
\end{gather}
Note that $+_{e}$ is not in general an associative operation, since
\begin{gather*}
\left(\left[-,e,b\right]+_{e}\left[-,e,d\right]\right)+_{e}\left[-,e,f\right]=\left[-,e,\left[\left[b,e,d\right],e,f\right]\right];\\
\left[-,e,b\right]+_{e}\left(\left[-,e,d\right]+_{e}\left[-,e,f\right]\right)=\left[-,e,\left[b,e,\left[d,e,f\right]\right]\right],
\end{gather*}
but if $\kappa$ is associative then $+_{e}$ is obviously associative
as well.

If $\kappa$ is commutative then $\left[x,e,\left[b,e,d\right]\right]=\left[x,e,\left[d,e,b\right]\right]$,
so 
\[
\left[-,e,b\right]+_{e}\left[-,e,d\right]=\left[-,e,d\right]+_{e}\left[-,e,b\right],
\]
 so the binary operation $+_{e}$ is commutative as well.

As $\left[x,e,\left[b,e,e\right]\right]=\left[x,e,b\right]$ by (\ref{eq:xyy})
and $\left[x,e,\left[e,e,b\right]\right]=\left[x,e,b\right]$ by (\ref{eq:xxy}),
we have 
\begin{equation}
\left[-,e,b\right]+_{e}\left[-,e,e\right]=\left[-,e,e\right]+_{e}\left[-,e,b\right]=\left[-,e,b\right]\label{eq:pointid}
\end{equation}
 for any $e,b\in X$, so $\left[-,e,e\right]$ is the unique identity
under pointed addition of translations.

Also, if (\ref{eq:k3}) holds then$\left[x,e,\left[b,e,\left[e,b,e\right]\right]\right]=\left[x,e,e\right]$,
so $\left[-,e,\left[e,b,e\right]\right]$ is a right inverse of $\left[-,e,b\right]$,
and if (\ref{eq:k4}) holds then $\left[x,e,\left[\left[e,b,e\right],e,b\right]\right]=\left[x,e,e\right]$,
so $\left[-,e,\left[e,b,e\right]\right]$ is a left inverse of $\left[-,e,b\right]$.
Thus, if (\ref{eq:k3}) and (\ref{eq:k4}) hold then
\begin{equation}
\left[-,e,b\right]+_{e}\left[-,e,\left[e,b,e\right]\right]=\left[-,e,\left[e,b,e\right]\right]+_{e}\left[-,e,b\right]=\left[-,e,e\right].\label{eq:pointinv}
\end{equation}
In particular, if $\kappa$ is associative then $\left[-,e,\left[e,b,e\right]\right]$
is the unique inverse of $\left[-,e,b\right]$ under pointed addition
of translations relative to $e$.

Recall that if (\ref{eq:xxy}) and (\ref{eq:pxxyz}) hold then $\left[-,a,b\right]=\left[-,e,b'\right]$
if and only if $b'=\left[e,a,b\right]$. This means that every translation
$\left[-,a,b\right]$ is a pointed translation which can be written
in the form $\left[-,e,b'\right]$ in exactly one way, namely as $\left[-,e,\left[e,a,b\right]\right]$.
(\ref{eq:defaddtl2}) can then be written as
\begin{gather*}
\left[-,a,b\right]+_{e}\left[-,c,d\right]:X\rightarrow X,\\
\left(x\mapsto\left[x,e,\left[\left[e,a,b\right],e,\left[e,c,d\right]\right]\right]\right)=\left[-,e,\left[\left[e,a,b\right],e,\left[e,c,d\right]\right]\right].
\end{gather*}
This implies that (\ref{eq:pointinv}) can be rendered as
\[
\left[-,a,b\right]+_{e}\left[-,b,a\right]=\left[-,b,a\right]+_{e}\left[-,a,b\right]=\left[-,e,e\right].
\]

Note that if (\ref{eq:pxxyz}) and (\ref{eq:pxyyz}) hold then 
\begin{gather*}
\left(\left[-,a,b\right]+_{e}\left[-,c,d\right]\right)\left(x\right)=\left(\left[-,a,b\right]+_{e}\left[-,b,\left[b,c,d\right]\right]\right)\left(x\right)\\
=\left[x,e,\left[\left[e,a,b\right],e,\left[e,b,\left[b,c,d\right]\right]\right]\right]=\left[x,e,\left[\left[e,a,b\right],b,\left[b,c,d\right]\right]\right]\\
=\left[x,e,\left[e,a,\left[b,c,d\right]\right]\right]=\left[x,a,\left[b,c,d\right]\right]=\left[\left[x,a,b\right],b,\left[b,c,d\right]\right]\\
=\left[\left[x,a,b\right],c,d\right]=\left(\left[-,a,b\right]+\left[-,c,d\right]\right)\left(x\right).
\end{gather*}
 That is, if $\kappa$ is an associative Malcev operation then 
\[
\left[-,a,b\right]+_{e}\left[-,c,d\right]=\left[-,a,b\right]+\left[-,c,d\right]
\]
 for any $e\in X$.\smallskip{}

\paragraph*{A5}

There is another, more well-known way of recovering a binary operation
from $\kappa$ satisfying (\ref{eq:xyy})--(\ref{eq:pxyyz}), this
time primarily involving individual points. For any fixed $e\in X$
the \emph{pointed combination} of $x$ and $y$ is the point $x\diamond_{e}y$
defined by 
\[
X\times X\rightarrow X,\qquad\left(x,y\right)\mapsto x\diamond_{e}y:=\left[x,e,y\right].
\]
If $\kappa$ is associative then $\left[\left[x,e,y\right],e,z\right]=\left[x,e,\left[y,e,z\right]\right]$,
so $\left(x\diamond_{e}y\right)\diamond_{e}z=x\diamond_{e}\left(y\diamond_{e}z\right)$. 

In general, $\left[x,e,e\right]=x$ for any $x\in X$ by (\ref{eq:xyy})
so that $e$ is a right identity element, and $\left[e,e,x\right]=x$
for any $x\in X$ by (\ref{eq:xxy}) so that $e$ is a left identity
element. Hence, if $\kappa$ is a Malcev operation then $e$ is the
unique identity in $\left(X,\diamond_{e}\right)$. Also, if (\ref{eq:k3})
holds then $\left[x,e,\left[e,x,e\right]\right]=e$ so that $\left[e,x,e\right]$
is a right inverse of $x$, and if (\ref{eq:k4}), holds then$\left[\left[e,x,e\right],e,x\right]=e$,
so that $\left[e,x,e\right]$ is a left inverse of $x$. In particular,
if $\kappa$ is associative then $\left[e,x,e\right]$ is the unique
inverse of $x$. 

This means that if $\kappa$ is an associative Malcev operation then
$\left(X,\diamond_{e}\right)$ is a group. If $\kappa$ is commutative
then $\left(X,\diamond_{e}\right)$ is evidently commutative.

Note that for any choice of identity element $e$ in $X$ there is
a bijection
\[
\phi:X\rightarrow\mathscr{T}_{X},\qquad p\mapsto\left[-,e,p\right].
\]
Recall that $\left[-,e,p\right]+_{e}\left[-,e,q\right]=\left[-,e,\left[p,e,q\right]\right]$.
Thus, $\phi\!\left(p\diamond_{e}q\right)=\left[-,e,\left[p,e,q\right]\right]=\left[-,e,p\right]+_{e}\left[-,e,q\right]=\phi\!\left(p\right)+_{e}\phi\!\left(q\right)$. 

Also, $\phi\!\left(e\right)=\left[-,e,e\right]$ so that $\phi\!\left(e\right)+_{e}\left[-,e,p\right]=\left[-,e,\left[e,e,p\right]\right]=\left[-,e,p\right]$
by (\ref{eq:xxy}) and $\left[-,e,p\right]+_{e}\phi\!\left(e\right)=\left[-,e,\left[p,e,e\right]\right]=\left[-,e,p\right]$
by (\ref{eq:xyy}). 

If $\kappa$ is an associative Malcev operation then 
\begin{gather*}
\phi\left(p^{-1}\right)\left(x\right)=\left[x,e,\left[e,p,e\right]\right]=\left[x,p,e\right]=\left[-,p,e\right]\left(x\right);\\
\left[-,e,p\right]+_{e}\phi\left(p^{-1}\right)=\left[-,e,p\right]+\phi\left(p^{-1}\right)=\left[\left[-,e,p\right],p,e\right]=\left[-,e,e\right];\\
\phi\left(p^{-1}\right)+_{e}\left[-,e,p\right]=\phi\left(p^{-1}\right)+\left[-,e,p\right]=\left[\left[-,p,e\right],e,p\right]=\left[-,p,p\right]=\left[-,e,e\right].
\end{gather*}

Thus, if $\kappa$ is associative then $\phi$ is a canonical isomorphism
between $\left(X,\diamond_{e}\right)$ and $\left(\mathscr{T}_{X},+_{e}\right)$.

Furthermore, if $\kappa$ is associative there is a canonical isomorphism
\[
\psi:\left(X,\diamond_{e}\right)\rightarrow\left(X,\diamond_{e'}\right),\qquad x\mapsto x\diamond_{e}e'
\]
for any $e,e'\in X$. In fact,
\[
\psi\left(x\diamond_{e}y\right)=\left(x\diamond_{e}y\right)\diamond_{e}e'=\psi\left(x\right)\diamond_{e'}\psi\left(y\right)
\]
because 
\begin{gather*}
\left[\left[x,e,y\right],e,e'\right]=\left[x,e,\left[y,e,e'\right]\right]=\left[\left[x,e,e'\right],e',\left[y,e,e'\right]\right].
\end{gather*}

Thus, from any associative Malcev operation on $X$ we can recover
a group $\left(X,\diamond\right)$, unique up to a canonical isomorphism.
This group can be interpreted as a group acting on $X$ by a group
action of the form $p\mapsto\left(-,e,p\right)$.

\end{document}